\documentstyle[aps,pra,multirow,hhline,epsfig]{revtex}     

\font\san=cmssbx10 at 12pt
\font\teneufm=eufm10
\font\seveneufm=eufm7
\font\fiveeufm=eufm5
\newfam\eufmfam
\textfont\eufmfam=\teneufm
\scriptfont\eufmfam=\seveneufm
\scriptscriptfont\eufmfam=\fiveeufm

\def\Lsc{{\cal L}}

\def\Vsc{{\cal V}}
\def\Ksc{{\cal K}}

\def\pbld{{\bf p}}

\def\qbf{{\bf q}}
\def\rbf{{\bf r}}
\def\vbf{{\bf v}}

\def\OO{{\Omega}}

\def\Tr{{\rm Tr}}

\def\vep{{\varepsilon}}

\def\Ksc{{\cal K}}

\def\Lsc{{\cal L}}

\def\Pss{{\mbox{\san P}}}

\font\san=cmssbx10

\newcommand{\beq}{\begin{equation}}
\newcommand{\eeq}{\end{equation}}
\newcommand{\beqa}{\begin{eqnarray}}
\newcommand{\eeqa}{\end{eqnarray}}

\begin{document} 

\begin{titlepage}

\title{Shock capturing by anisotropic diffusion 
        oscillation reduction}

\author{G. W. Wei}
\address{Department of Computational Science, 
National University of Singapore\\
Singapore 117543, R. Singapore}

\date{\today} 
\maketitle

\begin{abstract}
 
This paper introduces the method of anisotropic 
diffusion oscillation 
reduction (ADOR) for shock wave computations. 
The connection is made between digital image processing,
in particular, image edge detection,
and numerical shock capturing.
Indeed, numerical shock capturing can be 
formulated on the lines of iterative digital edge detection.
Various anisotropic diffusion  
and super diffusion operators originated from image
edge detection are proposed for 
the treatment of hyperbolic conservation laws
and near-hyperbolic hydrodynamic 
equations of change.  
The similarity between anisotropic diffusion and
artificial viscosity is discussed. 
 Physical origins and mathematical 
properties of the artificial viscosity is 
analyzed from the kinetic theory point of view.
A form of pressure tensor is derived from the 
first principles of the quantum mechanics. 
Quantum kinetic theory is utilized 
to arrive at macroscopic transport equations from the 
microscopic theory.  Macroscopic symmetry is used
to simplify  pressure tensor expressions. 
The latter provides a basis for the design of artificial 
viscosity. The ADOR approach is validated by using 
(inviscid) Burgers' equation in one and two spatial 
dimensions, the incompressible Navier-Stokes equation 
and the Euler equation. A discrete singular convolution (DSC) 
algorithm is utilized for the spatial discretization.

\end{abstract} 

\end{titlepage}

\section{Introduction}

Shock wave is a common phenomenon in nature, such as in 
aerodynamics and hydrodynamics. Mathematically, 
nonlinear hyperbolic conservation equations provide 
a good description to shock waves.
The construction of numerical schemes that are 
capable of shock capturing for hyperbolic and 
inviscid hydrodynamic equations is a major objective of
computational methodology. 
However, it is noted that 
the concept of discontinuity does not apply 
in digital computations.
Therefore, a shock in computational sense
may refer to systematic, rapid   
variation of function values over a few
grid points.   The difficulty is that 
hyperbolic equation may have weak solutions that are
discontinuous at the so-called shock front. 
Such a discontinuity will cause Gibbs' oscillations
in a high (spatial) order numerical scheme. The 
numerically induced oscillations are usually amplified
in (time) iterations. A variety of numerical schemes
have been proposed for shock capturing. As early as 50 
years ago, a solution to this problem was constructed 
by von Neumann and Richtmyer\cite{NeuRic}. The essence of 
their approach is to introduce  small artificial
viscosity so that a smooth solution can alway  be attained
in a finite difference approximation. A variety of 
modifications
to von Neumann and Richtmyer's method are 
made in the last 5 decades to address problems
of possible failure in spatial scaling and 
errors due to additional momentum flux and production,
as well as unbalanced heat flux and its over production in the 
simulation of hydrodynamic conservation laws.

A different approach was proposed by Godunov\cite{Godunov} 
to construct a full solution by using low order 
piecewise discontinuous approximations. 
Such a piecewise solution is a good approximation 
at the smooth regions, and is capable of representing
the shock front over a small region of grid.
The knowledge of  wave propagation and wave interaction is 
built in the numerical scheme in the form of a 
Riemann solver. The Godunov's approach has been 
extended to higher-order schemes\cite{Leer,WooCol81,WooCol84}.
In doing so, a higher-order approximation is constructed 
in the smooth region, while near the shock the solution 
is still of first order accuracy.  The Godunov method is very stable 
and thus, easy to design and use\cite{BelCol}.   
However the major disadvantage of this 
method is the complexity introduced into a numerical scheme 
through a Riemann solver. Such a drawback reduces
its computational efficiency. 

Another general approach is the hybrid scheme, 
which utilizes a high-order scheme for a smooth 
region while using a low-order scheme near a 
discontinuity. A linear combination of these two types of 
schemes is then used at each interface 
using weight factors which may be 
nonlinear functions of the local flow field. DeBar
constructed a linear hybridization of  
first-order and  second-order difference schemes as early as
1968\cite{Debar}. 
Harten and Zwas\cite{HarZwa} devised another early 
linear hybridization scheme. Boris and Book reported 
a blending algorithm which yields sharp discontinuities 
without oscillations\cite{BorBoo}.   
A total variation diminishing (TVD) scheme was proposed by 
Harten\cite{Harten83} to control the spurious oscillations 
in the numerical solution by using the total variation as a 
measure. The TVD scheme typically degenerates to first-order accuracy 
at locations with smooth extrema and was later generalized 
to an essentially non-oscillatory (ENO) scheme\cite{HEOC,ShuOsh89,EShu}. 
The major idea of the ENO scheme is to suppress spurious oscillations 
near the shock or discontinuities, while maintaining a 
higher-order accuracy at smooth regions. This line of thinking was 
further polished recently in a weighted  essentially 
non-oscillatory (WENO) scheme\cite{LOC}. 
The WENO approach takes a linear combination of a number of high-order 
schemes of both central difference  and  up-wind type. 
The central difference type schemes have a larger weight factor at the 
smooth region while the up-wind schemes play a major role at the shock 
or discontinuity. In general, these approaches are quite expensive 
since checks are performed before making 
a decision at each grid point. 
In terms of accuracy, all existing methods are at best of  first 
order near the shock or discontinuity.

Perhaps modified artificial viscosity methods are the most 
popular approaches in practical computations.  
Unfortunately, the artificial 
viscosity smears shocks over three or more grid zones, 
which can lead to serious errors in the physical 
interpretation of the numerical results. 
Special care is required to ensure that 
the smeared numerical shock is consistent with the true
thickness of the shock in a practical problem under study. 
This difficulty has led to enormous and continuous effort 
at developing efficient and robust approaches. 
One approach is to locally
refine the computational grid\cite{Gropp}. An alternative 
approach is the use of an adaptive 
mesh\cite{Eggleton,GDM}.
Both methods are aimed at  matching between the physical 
shock and the numerically smoothed one with respect to the spatial
extension. However, a major drawback is their restriction for 
extremely small time step sizes as required by the Courant 
constraint. In many cases, they have to be formulated in  an
implicit scheme, which imposes an extra complexity
in practical implementation and an extra requirement in
computer memory.  Another problem is the unbalanced momentum flux and 
production, and additional  heat flux and heat production 
in hydrodynamic equations.

The addition of a viscous term to the inviscid hydrodynamic 
equations or the hyperbolic equation has a physical 
justification, that the true physical flow has no 
discontinuities. Therefore, mathematical 
model can be modified so as to reflect the true physics.
In their original work, von Neumann and Richtmyer\cite{NeuRic}
have introduced an artificial viscosity $q_{NR}$ of the form
\begin{equation}\label{qnr}
q_{NR}\sim \left({\bf \nabla \cdot v}\right)^2
\end{equation}
for the equation of motion. Their term is of second 
order in gradient of the velocity field ${\bf v}$ (The
velocity gradient is a second rank tensor and has the 
divergence of velocity as its component)
and will not be very sensitive to small gradients. 
Landshoff\cite{Landshoff} 
proposed an additional term that vanishes 
less rapidly for small gradients
\begin{equation}
q_L\sim  {\bf \nabla \cdot v}.
\end{equation}
A generic form that contains both $q_{NR}$ and $q_{L}$
(i.e., $\alpha{\bf \nabla \cdot v}
+\beta \left(\nabla {\bf \cdot v}\right)^2$) 
has been carefully studied by  many 
researchers\cite{Kuropatenko,Wilkins,Schulz,CSW}.
In particular, Caramana et al\cite{CSW} discussed 
a number of intuitive criteria for this form and 
the parameterization
of artificial viscosity. Special considerations 
are given to dissipativity, Galilean invariance,
self-similar motion invariance, wave front invariance 
and viscous force continuity. 
This generic form, as it was proposed for hydrodynamic             
conservation laws, has its advantages over other 
approaches for certain fluid mechanical computations.
However, it is very unstable for certain hyperbolic 
equations due to the Courant-Friedrich-Lewy (CFL) 
stability condition.
For example, it does not work well in resolving 
a sharp shock front for inviscid Burgers' equation 
with a smooth initial condition.  We believe that this 
failure is due to the high-order in the gradient of
the Neumann and Richtmyer form, which  is too 
sensitive to large gradients at a sharp shock front. 
Despite the fact that a number of functional forms 
of the artificial viscosity have been proposed, the 
procedure is still ad hoc. For a given form, 
parameter selection is quite tricky and often 
varies from problem to problem.  

The existence of so 
many different approaches for shock capturing  
 indicates both the importance and the difficulty of 
the problem. The objective of the present paper is to introduce 
an anisotropic diffusion oscillation reduction (ADOR) 
approach for shock wave computations. 
The method of  anisotropic diffusion in association with a 
partial differential equation was proposed by 
Perona and Malik\cite{PeMa} for digital image processing in 1990.
Since then, much research has been stimulated in 
image processing and applied mathematics
communities\cite{CLMC,NiSh,ROF,WhPi,YXTK,ToTa,Shah,Acton,Kichenassamy,TBAB,BSMH}. 
The method uses the heat equation, which 
has a solution of the Gaussian type, as a low pass filter 
to eliminate noise in an image, 
while it detects and preserves the edges 
of the image.
It has been shown that the Perona-Malik equation provides
a computational approach to image segmentation, 
noise removal, edge detection, and image enhancement.
In a recent work, a generalized Perona-Malik 
equation was proposed 
for image restoration and edge enhancement\cite{weiieee99}.
In fact, digital edge detection has much in common with 
numerical shock capturing.
In this work, we introduce a set of generalized anisotropic
diffusion operators and edge enhancing functionals 
for shock wave computations. 

This paper is organized as follows: Section II is 
devoted to kinetic theory analysis of the 
equation of change and artificial viscosity. 
The physical origin and mathematical 
properties of artificial viscosity are
discussed from the kinetic theory point of view.
A set of hydrodynamic equations describing  fluid flow
consisting of microscopic particles is derived from the 
quantum Boltzmann equation, i.e., the Waldmann-Snider 
equation\cite{Waldmann,Snider60}. 
The latter can be regarded as a consequence of 
a reduction of the von Neumann equation, 
or the quantum Liouville equation
to the level of a single particle density operator.  
The mathematical form of the pressure tensor is analyzed
by using the group theory consideration. A set of generalized 
artificial viscosities, artificial heat and artificial 
virial correction are proposed as the results of the kinetic 
theory analysis. 
Generalized Perona-Malik equation is discussed in  
Section III for shock wave computations. This approach
follows a very different line of thinking as it was originally 
proposed for image analysis and processing.
An edge-controlled image enhancing functional is introduced 
for shock representation.
Anisotropic diffusion and diffusion super 
operators are introduced for 
fast and effective shock capturing.
Numerical experiments are presented in Section IV.
A few standard test problems, including 
Burgers' equation and 
inviscid Burgers' equation in one- and two-dimensions, 
the incompressible Euler and Navier-Stokes equations 
are used for exploring usefulness, testing efficiency
and illustrating the validity of the ADOR approach. 
A recently developed discrete singular convolution (DSC)
algorithm\cite{weijcp99,weiphysica20,weijpb20,weijpa20} 
is utilized for the numerical integration. 
This paper ends with a conclusion.

\section{Kinetic theory analysis of artificial viscosity}

Artificial viscosity was original proposed for shock capturing 
in association with hydrodynamic equations in the fluid mechanics.
The choice of its form and parameter should be consistent with 
and motivated by the physical origin of hydrodynamic 
equations. Since quantum theory is 
the foundation for the modern fluid mechanics, a
microscopic analysis is presented for appropriate understanding 
from the kinetic theory. The latter, such as the Boltzmann equation, 
serves as a theoretical basis for hydrodynamic conservation laws.

\subsection{Microscopic Analysis}

Consider the flow of a quantum gas system consisting of 
total $N$ particles in a volume $S$.  Its behavior is 
governed by the Schr\"{o}dinger equation 
\beq\label{sch1}
i\hbar {\partial \Psi\over\partial t}=H^{(N)}\Psi,
\eeq
where $H^{(N)}$ is the self-adjoint Hamiltonian of the 
system and $\Psi$ is a vector in the Hilbert space
associated with the system. 
For the description of physical observable, we adopt 
the density operator $\rho^{(N)}$ whose time evolution 
is governed by the quantum Liouville equation
\beq\label{sch2}
i{\partial \rho^{(N)}\over\partial t}=\Lsc^{(N)}\rho^{(N)}
={1\over\hbar}\left[H^{(N)},\rho^{(N)}\right]_-
={1\over\hbar}(H^{(N)}\rho^{(N)}-\rho^{(N)}H^{(N)}).
\eeq
A physical observable $O^{(N)}$ is a Hermitian operator of the Hilbert
space and has the expectation value given by 
\beq\label{sch2-2}
<O^{(N)}>={\rm Tr}_{1,\cdots,N}O^{(N)}\rho^{(N)},
\eeq
where ${\rm Tr}_{1,\cdots,N}$ is the trace over all the states of the 
$N$ particle system and in particular,
\beq
<1^{(N)}>={\rm Tr}_{1,\cdots,N}1^{(N)}\rho^{(N)}=1
\eeq
gives the normalization of the total probability for 
finding all of the $N$ particles in the volume $S$.
Here $1^{(N)}$ is the identity operator of the system.
A standard form for the Hamiltonian is given by 
\beq\label{sch3}
H^{(N)}=K^{(N)}+V^{(N)}=\sum_i K_i+
\sum_{i<j} V_{ij},
\eeq
where $K_i$ is the kinetic energy (including internal 
state energy) operator of particle $i$ and $V_{ij}$ is 
the potential operator of particles $i$ and $j$.
The physical behavior of the $N$-particle system is far too 
complicated to compute by any means as a macroscopic gas
flow may consist of $10^{23}$ particles or more. 
Fortunately, for ordinary gases, it is sufficient to 
consider the state  of a typical particle, say particle 1
\beq\label{sch4}
\rho^{(1)}_1=N {\rm Tr}_{2,\cdots,N}\rho^{(N)}.
\eeq
In general, 
 the state of $n$ particles is defined 
\beq\label{sch4-2}
\rho^{(n)}_{1,\cdots,n}=N(N-1)\cdots(N-n+1){\rm Tr}_{n+1,\cdots,N}\rho^{(N)}.
\eeq
The time evolution of particle 1 is governed by 
\beq\label{sch5}
i{\partial \rho^{(1)}_1\over\partial t}=\Lsc^{(1)}\rho^{(1)}_1
+{\rm Tr}_{2},\Vsc_{12}^{(2)}\rho^{(2)}_{12}
={1\over\hbar}\left[H^{(1)},\rho^{(1)}_1\right]_-
+{1\over\hbar}\left[V_{12},\rho^{(2)}_{12}\right]_-
\eeq 
where $H^{(1)}=K_1$.
This is the first member of the BBGKY hierarchy\cite{BBGKY,Bogo} 
in the quantum form. The general form of the BBGKY 
hierarchy\cite{BBGKY,Bogo} 
is given by 
\beq\label{sch6}
i{\partial \rho^{(n)}_{1,\cdots,n}\over\partial t}
=\Lsc^{(n)}\rho^{(n)}_{1,\cdots,n}
+{\rm Tr}_{2},\Vsc^{(n+1)}_{1,\cdots,n+1}\rho^{(n+1)}_{1,\cdots,n+1},
\eeq 
where $\Vsc^{(n+1)}_{1,\cdots,n+1}$ is the potential superoperator
between particles ${1,\cdots,n}$ and particle ${n+1}$.
This set of equations is formal and exact, since no approximation has been 
made. However, to determine the time evolution of  $\rho^{(1)}_1$,
it is necessary to know the behavior of $\rho^{(2)}_{12}$, which 
is, in turn, determined by the second order BBGKY equation 
and the latter involves three-particle density operator $\rho^{(3)}_{123}$.

\subsection{Mesoscopic analysis}

          Kinetic theory attempts to explain macroscopic
properties in terms of microscopic properties of the atoms 
and/or molecules based on the classical or quantum mechanics. 
Typical macroscopic observations deal with $N (\sim 10^{23})$
particles over a volume much larger than  the size of the individual
molecules, and over a time period much longer 
compared to the time scale of
the individual molecular dynamics. An ab-initio description of the time
dependence of such  macroscopic observations from the quantum Liouville
equation is practically impossible not only because of
the conceptual difficulty with the meaning of measurement, 
but also due to the large number of degrees of freedom involved. 
The kinetic theory approach simplifies the
$N$-particle problem dramatically by looking at the behavior of one
typical particle under the influence of all other particles.

      The most successful kinetic theory has been based on the Boltzmann
equation\cite{Boltz,Waldmann,Grad,Snider60}.
Specifically this has been related to the hydrodynamic equations and has
given rise to molecular expressions for the transport
coefficients. Various attempts have been made
to understand the Boltzmann equation from the $N$-body Liouville
equation, or utilizing  the BBGKY hierarchy described in the 
last subsection. A first principle ``derivation''
of the Boltzmann equation from the Liouville equation has been achieved
by Bogoliubov\cite{Bogo}, and independently by Green\cite{Green}.
    The Wang Chang-Uhlenbeck-de Boer equation\cite{WCUB} was
the first attempt to generalize the quantum Boltzmann equation, 
and to account for the internal degrees of freedom of a polyatomic 
gas. It was not
generally realized until 1957 that polyatomic gases can only be treated
properly by a quantum mechanical method due to the fact that the internal
energy levels of such molecules are, in general,
degenerate. Waldmann\cite{Waldmann} and, independently, 
Snider\cite{Snider60} introduced a quantum kinetic equation, 
the Waldmann-Snider equation,
\beq\label{Bolt}
i{\partial \rho^{(1)}_{1}\over \partial t} =
\Lsc^{(1)}_{1}\rho^{(1)}_{1}
+\Tr_{2}\Vsc^{(2)}_{12}
\OO_{L_{12}}\rho^{(1)}_{1}\rho^{(1)}_{2}
\eeq
where $\OO_{L_{12}}$ is the M\o ller superoperator for 
quantum mechanical scattering
\beq\label{moller}
\OO_{L_{12}}=
\lim_{t\rightarrow-\infty}e^{i\Lsc^{(2)}_{12}t}e^{-i\Ksc^{(2)}_{12}t}.
\eeq
The  Waldmann-Snider equation 
is still the basis for the kinetic theory of quantum 
gas flow\cite{MBKK,Snider90}. 
Note that a dramatically simplified version of the Boltzmann
equation provides the (lattice) Boltzmann gas (LBG) approach for the 
computational fluid dynamics in the recent years.

\subsection{Macroscopic analysis}

     The most important properties of a fluid flow are  
physical observables, or physical measurements. For conservative
observables, the equations of change of physical observables 
give rise to conservation laws. The equations of change are 
important not only because they govern the time dependencies 
of the macroscopic quantities but also because they are needed 
for  solving kinetic equations using the Chapman-Enskog 
method\cite{ChaCow}. The  equations of change for
physical observables are derived in this subsection.

      A single-particle observable for particle 1
is denoted by $\psi_1$. The expectation value of the physical 
observable $\psi_1$ is
determined by the singlet density operator $\rho^{(1)}_1$
according to
\beq\label{psi1c}
\langle \psi_1 \rangle
=\Tr_1\psi_1 \rho^{(1)}_{1}.
\eeq
The most important physical observables are the number 
density $n_1$,
stream velocity $\vbf_1$ and kinetic energy $\vep^K_1$.
Firstly, the physical observable
for  the number density $n_1$
is  the Dirac delta distribution
$\psi_1=\delta_1\equiv\delta(\rbf-\rbf_1)$ so that
\beq
n_1 =\langle \delta_1 \rangle.
\eeq
Secondly, the stream velocity $\vbf_1$ is given by
\beq
\vbf_1 \equiv {1\over n_1}\Tr_1 {1\over 2m_1} (\pbld_1\delta_1
+\delta_1\pbld_1) \rho_{1},
\eeq
where  $m_1$ is the mass of particle 1.
Finally, the kinetic energy $\vep^K_1$ per particle is given by
\beq
\vep^K_1 \equiv{1\over n_1} \Tr_1\left( {1\over 8m_1}
(\pbld_1^2\delta_1+2\pbld_1\cdot\delta_1\pbld_1
+\delta_1\pbld_1^2) -{1\over 2}m_1\vbf_1^2\delta_1 \right)\rho_1.
\eeq
These expressions are appropriately symmetrized to account for 
the fact that position and momentum do not commute, while
a physical observable is a Hermitian operator.
Although this operator approach is used here for 
for the number density $n_1$, stream
velocity $\vbf_1$ and kinetic energy $\vep^K_1$,
phase space expressions can be easily obtained 
for these fluid dynamic quantities by using the 
standard Wigner representation.  
The equations of change are of course independent of 
the detailed method used for the expression of  
various quantities.
Equations of change for a particle observable is obtained from
the quantum Boltzmann equation
 (\ref{Bolt})  according to
\beqa\label{expect1}
\hspace*{-1cm}{\partial \langle \psi_1 \rangle_1 \over\partial t}
- \langle {\partial \psi_1  \over \partial t}\rangle_1
={1\over i\hbar}\langle [\psi_1, H^{(1)}]_{-}\rangle_1+ {1\over i} 
\Tr_{1,2}
\Vsc^{(2)}_{12}
\OO_{L_{12}}\rho^{(1)}_{1}\rho^{(1)}_{2}.
\eeqa
In case that a physical observable $\psi_1$ is time independent, the
second term on the left hand side of Eq. (\ref{expect1}) vanishes.
Conservation laws are the equation of continuity
\beq\label{eoc3}
{\partial n_{1}\over\partial t}=-\nabla\cdot(n_1\vbf_1);
\eeq
the equation of motion
\beq\label{eom3}
n_1m_1{\partial\vbf_1\over\partial t}+n_1m_1\vbf_1\cdot\nabla\vbf_1
=-\nabla\cdot\Pss_1,
\eeq
and the kinetic energy equation
\beq\label{eoK1}
{\partial n_1\vep^K_1\over\partial
t}=-\nabla\cdot(n_1\vbf_1\vep^K_1+\qbf^K_1
+\qbf_{\rm coll})-\Pss^t_1{\bf \colon}\nabla\vbf_1+\sigma^K,
\eeq
where
the superscript $t$ denotes the transpose
and the kinetic heat flux is
\beq
\qbf^K_1=\Tr_1\left\{\left({\pbld_1\over
  m_1}-\vbf_1\right){(\pbld_1-m_1\vbf_1)^2\over
  2m_1}\delta_1\right\}_s\rho_1 ,
\eeq
where $\{~~\}_s$ designates appropriately operator-symmetrized
quantities. 
Here the collisional heat flux is
\beqa
\qbf_{\rm coll}&\equiv&-{1\over 8m_1}\Tr_{12}\int_{-1}^1d\nu\left\{\delta
   \left[\rbf- {\rbf_1+\rbf_2 \over 2}-{\nu\over2}(
   \rbf_1-\rbf_2) \right]\right.\nonumber\\
&\times&\left.
   (\pbld_1+\pbld_2-2m_1\vbf_1)\cdot{\partial V_{12}\over
   \partial \rbf_1}\right\}_s
\OO_{L_{12}}\rho^{(1)}_{1}\rho^{(1)}_{2}
\eeqa
and the kinetic energy production  results from the
 collisional transfer of energy
\beqa
\sigma^K&=&-{1\over 8m_1}\Tr_{12}\left\{(\delta_1+\delta_2)
   (\pbld_1-\pbld_2)\cdot{\partial V_{12}\over \partial \rbf_1}\right.
   \nonumber\\&& \left.
  ~~~~~~ +
   {\partial V_{12}\over \partial\rbf_1}\cdot(\pbld_1-\pbld_2)
   \right\}_s
\OO_{L_{12}}\rho^{(1)}_{1}\rho^{(1)}_{2}.
\eeqa
The pressure tensor  
$\Pss^{K}_1$ has kinetic and collisional contributions
\beq
\Pss_1=\Pss^K_1+\Pss^{\rm coll}_1
\eeq
where the kinetic pressure tensor is 
\beq
\Pss^K_1\equiv{1\over 4m_1}\langle
\delta_1\pbld_1\pbld_1+\pbld_1\delta_1\pbld_1
+(\pbld_1\delta_1\pbld_1)^t+\pbld_1\pbld_1\delta_1\rangle_1-
n_1m_1\vbf_1\vbf_1
\eeq
and the  collisional pressure tensor is given by 
\beqa
\Pss^{\rm coll}_1&=&-{1\over 4}\Tr_{12}(\rbf_1-\rbf_2)\nabla_1V_{12}\nonumber\\
  &\times&\int_{-1}^{1}\delta\left({\rbf_1+\rbf_2 \over 2}+{\nu\over2}(
  \rbf_1-\rbf_2) -\rbf\right)d\nu 
\OO_{L_{12}}\rho^{(1)}_{1}\rho^{(1)}_{2}.
\eeqa
The transpose $\Pss^t_1$ of the pressure tensor $\Pss_1$ enters to
couple the convective energy to the internal kinetic energy per
particle $\vep^K_1$, heat flux contributions arise from kinetic
$\qbf^K_1$ and collisional $\qbf_{\rm coll}$ motion and finally there is
a production term $\sigma^K$ which involves the transfer between
kinetic and potential energies. Although the equations of continuity and
motion are consistent with the conservation of the number of particles
and total momentum in the absence of bound pairs, the
presence of the production term in the kinetic energy implies that the
total kinetic energy is not conserved because of the possible conversion
to potential energy due to the non-locality of the collisions (in the 
macroscopic sense).  Thus it is necessary to look at the equation 
of change for the potential energy per particle $\vep^V_1$
\beq
n_1\vep^V_1={1\over4}\Tr_{12}(\delta_1+\delta_2)V_{12}\rho^{(2)}_{12}.
\eeq
This is obtained in a manner consistent with pair
particle interactions and shown to be of the form
\beq\label{eoV}
{\partial n_1\vep^V_1\over\partial
t}=-\nabla\cdot(n_1\vbf_1\vep^V_1+\qbf^V_1)
-\sigma^K,
\eeq
where
\beq
\qbf^V_1={1\over 4m_1}\Tr_{12}\left[(\pbld_1-m_1\vbf_1)\delta_1
       +\delta_1(\pbld_1-m_1\vbf_1)\right]V_{12}
\OO_{L_{12}}\rho^{(1)}_{1}\rho^{(1)}_{2}.
\eeq
The production of potential and kinetic energy exactly cancels
and there is an added heat flux contribution $\qbf^V_1$ 
associated with the conductive potential energy flow.

\subsection{Pressure tensor and artificial viscosity}

The estimation of transport coefficients from the 
Boltzmann equation is quite complicated
and involves many aspects. The most important transport quantities 
for physical and engineering applications are 
the mass, momentum and energy. Sometimes angular momentum 
transport may also be important for certain physical phenomena
and thus an equation of change for 
angular momentum may be required. However, angular 
momentum conservation is rarely discussed in the context
of hydrodynamic conservation laws. 
This is because the angular momentum is not measured 
like other physical observables
in experiments. It is noted that all hydrodynamic 
conservation laws have the structure of convection, conduction
and production.
The equation of momentum conservation leads to 
shock wave as the divergence of the pressure 
tensor vanishes. The existence of shock wave devastates
the numerical simulation as noted by  
non Neumann and Richtmyer\cite{NeuRic}.
They introduced artificial viscosity of the form of 
$({\bf \nabla\cdot v})^2$ to overcome the numerical 
difficulty. In fact, a much better choice can be selected
by analyzing the form of the pressure tensor.

The Pressure tensor is a tensor of second rank and has 
two indices.
It has contributions from the kinetic transport 
$\mbox{\san P}^K$ 
and collision transport $\mbox{\san P}^{\rm coll}$ 
(Here, the subscript 1 is omitted for convenience).
For example, the kinetic transport $\mbox{\san P}^K_{xy}$
is the rate of transport (kg m/s)/(m$^2$s) in the x-direction,
of the momentum with respect to  the 
y-direction. The hydrodynamic pressure
\beq\label{p1}
\mbox{\san P}_{\rm eq}= nk_BT\mbox{\san U},
\eeq
given by the equation of state, is the 
local equilibrium contribution to the 
pressure tensor. Here, $\mbox{\san U}$ is the identity
tensor of  second rank and $k_B$ is the Boltzmann 
constant. In general, the gas is not at local equilibrium 
and the structure of the nonequilibrium part of the 
pressure tensor, 
\beq
{\bf \Pi}= \mbox{\san P}-\mbox{\san P}_{\rm eq}
\eeq
can be analyzed according to the irreducible representation 
of the improper three-dimensional rotational group $O(3)$.
As such,  ${\bf \Pi} $ can be expressed as a linear 
combination the irreducible  representations 
\beq
{\bf \Pi}=\Pi\mbox{\san U}+{\bf \varepsilon\cdot}\mbox{\san P}^a
+{\bf \Pi}^{(2)},
\eeq
where scalar, antisymmetric and symmetric traceless parts of 
the second rank tensor are given by 
\begin{eqnarray}\label{tensor}
\Pi&=&{1\over3}\mbox{\san U}{\bf :\Pi}\nonumber\\
\mbox{\san P}^a&=&-{1\over2}{\bf \varepsilon\cdot}{\bf \Pi}\nonumber\\
{\bf \Pi}^{(2)}&=&{1\over2}[{\bf \Pi}+{\bf \Pi}^t]
-{1\over3}\mbox{\san U}[\mbox{\san U}{\bf :} {\bf \Pi}].
\end{eqnarray}
Here ${\bf \varepsilon}$ is the Levi-Civita tensor,
$\Pi$ is a scalar (one dimension),
$\mbox{\san P}^a$ is a vector (three dimensions)
and ${\bf \Pi}^{(2)}$ is a tensor (five dimensions).
In the theory of the Boltzmann equation and 
irreversible thermodynamics,
the pressure tensor is treated as proportional to 
velocity gradient as in the case of 
a classical Newtonian flow. 
The velocity gradient, ${\bf \nabla v}$, 
is also a second rank tensor and can also be decomposed 
into three different components of
order zero, one and two, similar to Eq. (\ref{tensor}), 
\begin{eqnarray}\label{tensor1}
{\bf \nabla\cdot v}&=&{1\over3}\mbox{\san U}
{\bf :}\nabla{\bf v}\nonumber\\
{\bf \nabla\times v}&=&-{1\over2}{\bf \varepsilon\cdot}\nabla{\bf v}
\nonumber\\
\left[\nabla
{\bf v}\right]^{(2)}
&=&{1\over2}[\nabla{\bf v}+(\nabla{\bf v})^t]
-{1\over3}\mbox{\san U}\nabla\cdot{\bf v}.
\end{eqnarray}
Here, ${\bf \nabla v}$ is also a nine-dimensional quantity.
In principle, there are 81 possible coefficients connecting 
the pressure tensor and the velocity gradient. 
Symmetry analysis indicates that there
are three phenomenological equations
relating the corresponding components of ${\bf \Pi}$ and ${\bf \nabla v}$ 
\begin{eqnarray}\label{tensor2}
\Pi&=&-\eta^0_V{\bf \nabla\cdot v}\nonumber\\
\mbox{\san P}^a&=&-\eta^0_r{\bf \nabla\times v}\nonumber\\
{\bf \Pi}^{(2)}&=&-2\eta^0[{\bf \nabla v}]^{(2)},
\end{eqnarray}
where $\eta^0_V, \eta^0_r$ and $\eta^0$ are the 
viscosities of bulk, rotational and shear.

The pressure tensor given by Eq (\ref{tensor2}) provides 
an excellent approximation for a gas system near the 
local equilibrium, i.e. the Boltzmann distribution 
\beq
\rho^{(1)}_1={1\over {\rm Tr}_1e^{-K_1/k_BT}}e^{-K_1/k_BT}.
\eeq
For the purpose of numerical stability and shock wave 
simulation, we can impose artificial pressure tensor as 
\begin{eqnarray}\label{tensor3}
\Pi_{\rm art}&=&-\zeta^0_V{\bf \nabla\cdot v}\nonumber\\
\mbox{\san P}^a_{\rm art}&=&-\zeta^0_r{\bf \nabla\times v}\nonumber\\
{\bf \Pi}^{(2)}_{\rm art}&=&-2\zeta^0[{\bf \nabla v}]^{(2)},
\end{eqnarray}
where $\zeta^0_V, \zeta^0_r$ and $\zeta^0$ are artificial
bulk viscosity, artificial rotational viscosity 
and artificial shear viscosity, respectively.
These artificial viscosities modify the non equilibrium 
part of the pressure tensor. 
In fact, from the point of view of computations,
the equilibrium part of the 
pressure tensor can also be modified, i.e. a term which is 
proportional to $nk_BT$
\beq
\mbox{\san P}_{\rm nT}
=\epsilon nk_BT \mbox{\san U}.
\eeq
$\mbox{\san P}_{\rm nT}$ can be interpreted 
either as artificial heat or artificial virial correction.
Hence, we have the artificial pressure 
tensor of the form
\begin{eqnarray}\label{tensor3-2}
\mbox{\san P}_{\rm art} &=&\mbox{\san P}_{\rm nT}+
{\Pi}_{\rm art}\mbox{\san U}
+{\bf \varepsilon\cdot}\mbox{\san P}^a_{\rm art} 
+ {\bf \Pi}^{(2)}_{\rm art}\\
&=&
\epsilon nk_BT \mbox{\san U}
-\zeta^0_V{\bf \nabla\cdot v}\mbox{\san U}
-\zeta^0_r {\bf \varepsilon\cdot} {\bf \nabla\times v}
-2\zeta^0[{\bf \nabla v}]^{(2)}.
\end{eqnarray}
It should be  noted that there is an advantage in keeping
only a part of the artificial pressure tensor that 
corresponds to the non-equilibrium pressure tensor.
The latter makes no additional contribution 
to conservative quantities, such as, number density,
momentum and energy, as required by the Fredholm 
alternative for the existence of a solution for the 
kinetic equation in the Chapman-Enskog 
method\cite{ChaCow}.

However, for a gaseous system far from the local equilibrium,
higher-order terms in the velocity gradient become 
important. Large velocity gradient can build up from 
special boundary condition, such as the boundary layer 
phenomena, or from the lack of relaxation in a very
dilute gas flow. Therefore artificial viscosities and 
artificial heat of Eqs. (\ref{tensor3})
can be modified to include higher-order velocity gradient
contributions. Recognizing that the pressure tensor is 
a second rank tensor, thus all velocity gradients 
contribute in one of the forms of 
$\Pi, \mbox{\san P}^a$ and ${\bf \Pi}^{(2)}$. 
In this regard, 
the artificial viscosity form proposed by 
von Neumann and Richtmyer\cite{NeuRic}, 
$({\bf \nabla\cdot v})^2$, contributes to $\Pi$.
Certainly, artificial viscosities of the forms of 
${\bf \nabla\cdot v}$, ${\bf \nabla\times v}$ 
and $[{\bf \nabla v}]^{(2)}$ can be used for numerical 
computations. All of the aforementioned forms are 
at least of the first order in gradient, which may lead
a strong smoothing at the edge of a shock front 
and thus lead to large errors in numerical applications. 
Furthermore, these forms can result in instability when the 
shock front is very sharp, such as in  inviscid 
Burgers' equation. 
For these reasons, it is appropriate to 
propose  a general expression for the artificial 
viscosity and artificial heat so as to allow all
coefficients to be functions of  
${\bf \nabla\cdot v}$, ${\bf \nabla\times v}$ 
and $[{\bf \nabla v}]^{(2)}$ 
\begin{eqnarray}\label{tensor4}
\epsilon&=& \epsilon
\left({\bf \nabla\cdot v},{\bf \|\nabla\times v\|},\|[{\bf \nabla v}]^{(2)}\|\right)
\\ 
\zeta^0_V&=& \zeta^0_V
\left({\bf \nabla\cdot v},{\bf \|\nabla\times v\|},\|[{\bf \nabla v}]^{(2)}\|\right)
\\
\zeta^0_r&=& \zeta^0_r  
\left({\bf \nabla\cdot v},{\bf \|\nabla\times v\|},\|[{\bf \nabla v}]^{(2)}\|\right)
\\
\zeta^0 &=& \zeta^0  
\left({\bf \nabla\cdot v},{\bf \|\nabla\times v\|},\|[{\bf \nabla v}]^{(2)}\|\right)
\end{eqnarray}
where $\|\cdot \|$ denotes the magnitude and is computed as
\beq
{\bf \|A\|}=\sqrt{\sum_iA_i^2}
\eeq 
for a vector ${\bf A}$, and 
\beq
{\|\mbox{\san B}\|}=\sqrt{\sum_{ij}\mbox{\san B}_{ij}^2}
\eeq 
for a tensor $\mbox{\san B}$.
Obviously, artificial viscosity of 
von Neumann and Richtmyer\cite{NeuRic}, $({\bf \nabla\cdot v})^2$,
becomes a special case of the present treatment. 
The von Neumann and Richtmyer form can be classified either 
as an artificial heat term, $\mbox{\san P}_{\rm nT}$, 
or as an artificial bulk viscosity term,
${\Pi}_{\rm art}\mbox{\san U}$.

The inclusion of an artificial rotational viscosity can 
be an efficient way for the treatment of fast rotational  flow. 
It is noted that the rotational viscosity is in the 
original Navier-Stokes equation derived from
the kinetic theory. However, for flows with natural
boundary conditions, the angular momentum is also conserved
and it does not couple  the linear momentum except through
internal spinor relaxations. However, this is 
no longer the case in a flow with an irregular geometry.
For such a case,
the internal-conversion  between angular momentum and
linear momentum is substantial and therefore the rotational 
viscosity, $\eta^0_r$, should be retained in both 
theoretical modeling and numerical simulation. 
A detailed discussion and numerical
simulation of this aspect will be accounted elsewhere.
In the next section, we analyze the shock capturing 
algorithm further from the point of view of 
image processing, particularly, image edge detection.

\section{Oscillation reduction by anisotropic diffusion}

Although anisotropic diffusion was original proposed 
for image processing, 
there is much in common between digital 
image processing and computational fluid dynamics. 
An image function $I({\bf r})$ is   
a two-dimensional projection  of certain    
physical quantities, such as matter, velocity, energy,    
electromagnetic field, etc., under appropriate   
illumination conditions.
The edges in an image usually refer to rapid changes   
in some physical properties, such as geometry,   
illumination, and reflectance. Mathematically, a   
discontinuity may be involved in the function   
representing such  physical properties. 
Therefore image edges are very similar to shocks in 
fluid dynamics. Numerical shock capturing can be
formulated on the lines of iterative  
digital edge detection.
Edge detection is a key issue in  image processing,   
computer vision, and pattern recognition.   
A variety of algorithms, such as
the Sobel operator, the Prewitt operator,   
the Canny operator\cite{canny} and the 
DSC algorithm\cite{HouWei}
are proposed for image edge 
detection and representation. 
Anisotropic diffusion is a promising new mathematical 
algorithm  for image edge detection and image processing.  
This  basic idea can be adopted and modified for numerical 
shock capturing in association with the hyperbolic 
conservation laws. The Perona-Malik algorithm\cite{PeMa}
is reviewed
before the method of anisotropic diffusion oscillation 
reduction (ADOR) is discussed.

\subsection{The Perona-Malik equation}

The basic idea behind the Perona-Malik algorithm is to 
evolve an original image, $I({\bf r})$, under an 
edge-controlled diffusion operator\cite{PeMa}
\begin{eqnarray}\label{PM}
\frac{\partial u({\bf r},t)}{\partial t}
&=&{\bf \nabla}\cdot \left[d(\|\nabla u({\bf r},t)\|)
\nabla u({\bf r},t)\right]\nonumber\\
u({\bf r},0)&=&I({\bf r}).
\end{eqnarray}
Here, $d(\|\nabla u\|)$ is a generalized diffusion coefficient
which is so designed that its values are very small at the edges 
of an image. Many edge stopping functions $d(|\nabla u|)$ 
are appropriate for anisotropic diffusion. For example, 
the Gaussian 
\beq
d(\|\nabla u\|)=e^{-|\nabla u|^2/2\sigma^2} 
\eeq
and the Lorentz
\beq
d(\|\nabla u\|)={1\over 1+|\nabla u|^2/\sigma^2} 
\eeq
are both suitable for edge representation. They  provide
perceptually similar results in practical applications. 
Numerically, the diffusion 
coefficient becomes very small near an image edge due to 
the effect of edge stopping functions $d(\|\nabla u\|)$.
As a result, the image edge is preserved in the diffusion 
process. The pixel values at a non-edge  part will be 
smoothed and reduced due to the substantial 
diffusion coefficient
prescribed by $d(\|\nabla u\|)$.
Perona and Malik argued that the solution of 
their anisotropic 
diffusion equation  has no additional maxima 
(minima) which does not belong to the 
initial image data. However, this point has been challenged 
recently\cite{CLMC,Kichenassamy}.
It is well-known
that this anisotropic diffusion algorithm may 
break down when the gradient generated by noise is 
comparable to image edges and features. 
Numerically, this can be alleviated by using 
a regularization procedure.

\subsection{Shock capturing by anisotropic diffusion} 

   Hyperbolic conservation laws of the type
\beq\label{hyper}
\frac{\partial u({\bf r},t)}{\partial t}
+{\bf \nabla\cdot F}(u)=0
\eeq  
describe the rate of change  of a physical quantity
$u$ given by the generalized convection 
${\bf \nabla\cdot F}(u)$. 
Without the balance of conduction and/or production, 
Eq. (\ref{hyper}) may have a discontinuous solution.
The task is to construct a stable computational 
scheme which is capable of resolving the ``shock''. 
We first note that computationally,  
if a shock is defined as a 
discontinuity, there is no shock to capture. 
This is because,  
the original notion of discontinuity is 
undefined on a discrete mesh. Therefore, 
a numerical shock is characterized by rapid 
variation of function values over a small 
grid zone. Numerical shock capturing, at best, is globally 
a first order approximation scheme for resolving 
the large gradient feature of the true solution.
However, locally, it is preferred
to compute the solution as accurate as possible, 
except at the discontinuity of the true solution.
To this end, we introduce an anisotropic diffusion 
term to Eq. (\ref{hyper})
\begin{eqnarray}\label{hyper2}
\frac{\partial u({\bf r},t)}{\partial t}
+{\bf \nabla\cdot F}(u)&=&
{\bf \nabla}\cdot 
\left[d_1(\|\nabla u({\bf r},t)\|)\nabla u({\bf r},t)\right]
\end{eqnarray}
where the diffusivity, $d_1(\|\nabla u({\bf r},t)\|$, 
is chosen such that it is essentially zero except at 
a numerical shock position.
Obviously, $\|\nabla u({\bf r},t)\|$
is important for shock detection.
Apparently, Eq. (\ref{hyper2}) reduces to the 
Perona-Malik equation (\ref{hyper}) without the convection term.
However, the selections of the anisotropic diffusivity in these two
equations are entirely different and
 serve   opposite purposes. 
For example, one may choose $d_1(\|\nabla u({\bf r},t)\|$
as
\beq\label{c1-0}
d_1=d^1_1 \ln[(\nabla u)^2+1],
\eeq
where $d^1_1$ is a constant. Equation (\ref{c1-0})
differs much from the Gaussian and the Lorentz.

The derivation of Eq. (\ref{hyper2}) has its roots in the
conservation of a physical quantity involving a
phenomenological flux and its divergence. 
However, as a computational algorithm,  
an anisotropic diffusion of the form  
\begin{eqnarray}\label{hyper3}
\frac{\partial u({\bf r},t)}{\partial t}
+{\bf \nabla\cdot F}(u)&=&
\Gamma_1(\|\nabla u({\bf r},t)\|)
\nabla^2 u({\bf r},t)
\end{eqnarray}
can be used. In fact, expressions in both Eqs. (\ref{hyper2}) 
and (\ref{hyper3}) are efficient for numerical shock capturing.
Here, $\Gamma_1$ is chosen to smooth the  
oscillations near a shock and is essentially zero in other
regions.

\subsection{Edge enhancing functional and super diffusion}

From the hydrodynamic point of view, the governing equation of
a conservative quantity has a general structure, i.e., the
rate of change is balanced by convection,
conduction and production. Therefore, the 
nonlinear advective motion can be counterbalanced 
by an appropriate production. Therefore, 
we propose a real-valued, bounded shock (edge)
enhancing functional  
\begin{equation}\label{enh}
e(\|\nabla u({\bf r},t)\|).
\end{equation}
Here $e(\|\nabla u({\bf r},t)\|)$ 
is appropriately 
chosen so that it is edge sensitive and is essentially 
zero away from a numerical shock or ``discontinuity''.
This leads to another shock capturing equation
\begin{eqnarray}\label{GPM1}
\frac{\partial u({\bf r},t)}{\partial t}
+{\bf \nabla\cdot F}(u)
&=&{\bf \nabla}\cdot \left[d_1(\|\nabla u({\bf r},t)\|)
\nabla u({\bf r},t)\right]\nonumber\\
&+&
e(\|\nabla u({\bf r},t)\|). 
\end{eqnarray}         
Appropriate choice of the shock enhancing functional
will result in a stable numerical algorithm. 
Obviously, the von Neumann and Richtmyer\cite{NeuRic}
artificial viscosity term, Eq. (\ref{qnr}), 
is a special case of 
the present formulation.

   The diffusion equation can be derived from 
Fick's law for mass flux, 
\beq
{\bf j_1}({\bf r},t)=-D_1 \nabla u({\bf r},t)
\eeq
with $D_1$ being a constant. 
From the point of view of kinetic theory, 
this is an approximation to a quasi homogeneous 
system which is near equilibrium. A better approximation 
can be expressed as a super flux  
\begin{equation}\label{heat2}
{\bf j_q}({\bf r},t)=-\sum_q D_q\nabla \nabla^{2q}u({\bf r},t), 
~~~(q=1,2,\cdots),
\end{equation}
where $D_q$  are constants and 
higher order terms ($q>1$) describe corrections to mass
flux by the influence of inhomogeneity in density 
distribution and of flux-flux correlations.
The mass conservation leads to 
\begin{eqnarray}\label{heat3}
\frac{\partial u({\bf r},t)}{\partial t}&=&-\nabla\cdot
    {\bf j_q}({\bf r},t)+s({\bf r},t)\nonumber\\
&=&\sum_q \nabla\cdot[D_q\nabla\nabla^{2q}u({\bf r},t)]+s({\bf r},t),
~~~(q=1,2,\cdots),
\end{eqnarray}
where $s$ is a source term which can be a nonlinear 
function describing chemical reactions. 
Equation (\ref{heat3}) is a generalized reaction-diffusion 
equation which includes not only the usual diffusion  
and production terms, but also super diffusion terms. 
A truncation at the second order super flux, 
${\bf j_2}({\bf r},t)=-D_1\nabla u({\bf r},t)
-D_2\nabla\nabla^{2}u({\bf r},t)$, leads to the expression
that is important in many phenomenological theories, such as the
Cahn-Hilliard equation and the Kuramoto-Sivashinsky 
equation. The latter have been used for the 
description of a number of physical phenomena, such as the
Taylor-Couette flow, parametric waves and pattern formation 
in alloys, glasses, polymers, combustion and biological 
systems\cite{CroHoh}. 

In a shock wave simulation process, the distribution 
of the physical quantity under study can be 
highly inhomogeneous and/or oscillatory. Hence, the present
 shock capturing algorithm can be made
more efficient by incorporating a shock (edge) sensitive 
super diffusion operator
\begin{eqnarray}\label{GPM2}
\frac{\partial u({\bf r},t)}{\partial t}
+{\bf \nabla\cdot F}(u)
&=&\sum_q{\bf \nabla}\cdot \left[d_q(\|\nabla u({\bf r},t)\|)
\nabla\nabla^{2q} u({\bf r},t)\right]\nonumber\\
&+& e(\|\nabla u({\bf r},t)\|),
~~~(q=1,2,\cdots).
\end{eqnarray}
Here $d_q(u,|\nabla u|)$ are shock (edge) sensitive diffusion functions. 
In most applications, a truncation at $q=2$ is sufficient  
\begin{eqnarray}\label{GPM3}
\frac{\partial u({\bf r},t)}{\partial t}
+{\bf \nabla\cdot F}(u)
&=&{\bf \nabla}\cdot [d_1(\|\nabla u({\bf r},t)\|)
\nabla u({\bf r},t)]
\nonumber\\
&+&{\bf \nabla}\cdot [d_2(\|\nabla u({\bf r},t)\|)\nabla\nabla^{2} 
u({\bf r},t)]
\nonumber\\
&+&e(\|\nabla u({\bf r},t)\|).
\end{eqnarray}
This equation can be simplified further 
for practical implementation
\begin{eqnarray}\label{GPM4}
\frac{\partial u({\bf r},t)}{\partial t}
+{\bf \nabla\cdot F}(u)
&=&
\Gamma_1(\|\nabla u({\bf r},t)\|)\nabla^2 u({\bf r},t)
\nonumber\\
&+&
\Gamma_2(\|\nabla u({\bf r},t)\|)\nabla^{4} 
u({\bf r},t)
\nonumber\\
&+&e(\|\nabla u({\bf r},t)\|),
\end{eqnarray}
where $\Gamma_1$ and $\Gamma_2$ are to be appropriately chosen
to capture the shock. Usually, $\Gamma_1$ is a positive 
function, while $\Gamma_2$ is negative.

From the point view of image processing, both  operators 
$\nabla^{2}$ and  $\nabla^{4}$ are high pass filters.
In the Fourier representation, operators 
$\nabla^{2}$ and  $\nabla^{4}$
are proportional to $\omega^2$ and $\omega^4$ of the 
frequency $\omega$. 
Hence, $\nabla^{2}$ is more sensitive to low frequency
oscillations while $\nabla^{4}$ has a large filter 
response for high frequency oscillations. Although both 
$\nabla^{2}$ and  $\nabla^{4}$ become band pass 
filters on a discrete grid, their filter responses have 
similar properties as in their abstract operator forms.

Note that many of these anisotropic terms presented 
in this section resemble those expressions
derived in the last section for the artificial pressure
tensor, $\mbox{\san P}_{\rm art}$. Although the 
starting points for these  two approaches are entirely 
different, the results are obviously connected,
notation needs to be simplified and
possible confusions are to be avoided. 
For these reasons, we use the acronym 
ADOR for both approaches.

\section{Numerical experiments}       

In this section, the numerical schemes developed 
from the analysis of kinetic theory and  image processing 
analysis are utilized for numerical computations 
involving shock waves. Parameter selections for 
the ADOR method are specified. Obviously,
there are infinitely many ways to construct anisotropic 
diffusion functionals. For simplicity, we test the 
following choices of the ADOR coefficients: 
\beq\label{c1}
d_1=d^1_1 \ln[(\nabla u)^2+1],~ d_2=d_2^1\ln[(\nabla u)^2+1],
~e=e^1(\nabla\cdot u)^2 
\eeq
for Eq. (\ref{GPM3}),
\beq\label{c2}
\Gamma_1=\gamma^1_1(|u_x|)^{1\over4},~
\Gamma_2=\gamma^1_2(|u_x|)^{1\over4},~
e=0
\eeq
for Eq. (\ref{GPM4}),
\beq\label{c3}
\Gamma_1=\gamma^2_1\ln[|u_x|+1],~
\Gamma_2=\gamma^2_2\ln[|u_x|+1],~     
~e=0, 
\eeq
for Eq. (\ref{GPM4}), 
\beq\label{c4}
\Gamma_1=\gamma^3_1\ln[u^2_x+1],~
\Gamma_2=\gamma^3_2\ln[u^2_x+1],~     
~e=0, 
\eeq
for Eq. (\ref{GPM4}), and 
\beq\label{c5}
\Gamma_1=\gamma^4_1\ln[\|\nabla u\| +1],~
\Gamma_2=\gamma^4_2\ln[\|\nabla u\|  +1],~     
~e=0, 
\eeq
for Eq. (\ref{GPM4}).
In fact, coefficients in
x- and y-direction can be chosen separately
\beq\label{c6}
\Gamma_1= \gamma^5_1\left\{\begin{array}{ll}
   \ln[|u_x|+1]~~ &\mbox{for $u_{xx}$}\\
   \ln[|u_y|+1]~~ &\mbox{for $u_{yy}$}\\
   \end{array}\right., 
~\gamma^5_2=0,~e=0. 
\eeq 
Note that these expressions have one feature in common, i.e.,
they are all of low order in the gradient of $u$. This feature
allows a sharp change in the solution to be handled by 
a coarse grid. Otherwise, if there are higher-order
gradient terms,  such as the 
von Neumann and Richtmyer form\cite{NeuRic},
 the computational grid near the shock front has
to be refined to reduce the flux amplitude. A finer 
grid, in turn, leads to another stability problem as 
dictated by the CFL condition.

In the rest of this section,
the performance of the abovementioned prescriptions 
is examined for a detailed choice of ADOR parameters
$d^1_1,d^1_2,e^1$ and $\gamma^i_j ~(i=1,2,3;j=1,2)$.
A few standard problems are employed to
test the validity, and to demonstrate 
the robustness of the present approach.
These problems include Burgers' equation 
in one and two space dimension, the 
incompressible Navier-Stokes 
equation and Euler equation with periodic 
boundary conditions. 

For spatial discretization, we use 
the discrete singular convolution (DSC)    
algorithm which was proposed as a potential approach     
for computer realization of singular     
convolutions\cite{weijcp99,weiphysica20,weijpb20,weijpa20}. 
Mathematical foundation     
of the algorithm is the theory of   
distributions\cite{Schwartz}.   Sequence of     
approximations to the singular kernels of Hilbert     
type, Abel type and delta type were constructed.   
Applications are discussed to analytical    
signal processing, Radon transform and surface   
interpolation.        
Numerical solutions to differential equations     
are formulated via singular kernels of delta type.    
By appropriately choosing the DSC kernels, the     
DSC approach exhibits global methods' accuracy     
for integration and local methods' flexibility     
for handling complex geometries and boundary conditions. 
Unified features of the DSC approach to differential 
equations were analyzed in details.
In particular, we demonstrated\cite{weijpa20} 
that different implementations of the 
DSC algorithm, such as global, local,
Galerkin, collocation, and finite difference, 
can be deduced from a single starting point.    
The DSC algorithm is validated for
the numerical solution of     
the Fokker-Planck equation\cite{weijcp99},
the Schr\"{o}dinger equation\cite{weijpb20} and  
the Navier-Stokes equation\cite{weicmame20,Wanweip2}.
It was also utilized to integrate the     
(nonlinear) sine-Gordon equation     
with the initial values close to a homoclinic     
orbit singularity\cite{weiphysica20}, for which     
conventional local methods    
encounter great difficulties and   
numerically induced chaos was     
reported for such an integration\cite{Ablowitz}.
The reader  is referred to  \cite{weijcp99} for a 
detailed discussion of the method.   
The DSC parameters used in the present work are
$\sigma/\Delta=3.2$ and $M=31$ for the DSC kernel 
of regularized Shannon\cite{weijcp99}. 
For time discretization, the explicit
4th-order Runge-Kutta scheme is used for 
Burgers' equation. The Navier-Stokes equation
is integrated by using the implicit Euler scheme
in association with a discrete singular 
convolution-alternating direction implicit 
(DSC-ADI) algorithm\cite{zhaoADI}. 
Further information is given in 
subsections below.

\subsection{Burgers' equation in one dimension}

Burgers' equation\cite{Burger} is an important model 
for the understanding of physical flows.
It appears customary to test 
new schemes in computational fluid dynamics 
by applying them to Burgers' equation. 
Despite  much of the effort, numerical 
solution  of  Burgers' equation is still not a trivial 
task, particularly at very high Reynolds numbers
where the nonlinear advection leads to shock waves.
In fact, many standard
computational algorithms fail to predict 
Burgers' inviscid shocks.

Burgers' equation is  given by 
\begin{equation}\label{Burgers}
{\partial u\over \partial t}+u{\partial u\over\partial x}
={1\over {\rm Re}}  
{\partial^2 u\over\partial x^2}, 
\end{equation}                 
where $u(x,t)$ is the dependent variable 
resembling the flow velocity and Re is the 
Reynolds number characterizing the size of the viscosity.
The competition between the nonlinear advection and the
viscous diffusion is 
controlled by the value of Re in Burgers' equation, and
thus determines the behavior of the solution. 
 We consider Eq. (\ref{Burgers}) using the 
following initial and boundary conditions
\begin{eqnarray}\label{bur2}
u(x,0)&=&\sin(\pi x),\nonumber\\
u(0,t)&=&u(1,t)=0.
\end{eqnarray}
Cole has provided an exact solution\cite{Cole} for this problem 
in terms of a series expansion which is readily 
computable roughly for the parameter ${\rm Re}\leq 100$.     
For the parameter ${\rm Re}=100$, the present calculations use 
41 grid points in the interval [0,1]
with a time increment of 0.01. 
Both $L_{1}$  and $L_{\infty}$ errors
at 11 different times are listed in TABLE I.

We next consider inviscid Burgers' equation
with same initial and boundary condition as given in 
Eq. (\ref{bur2}). In this case, the solution quickly 
develops into a sharp shock front at $x=1$. 
The DSC algorithm does not work along because severe
oscillations eventually turn into an uncontrolled 
error growth. This problem is treated by using 
the anisotropic diffusion approach discussed 
in the previous sections. The performance of 
four different prescriptions given in 
Eqs. (\ref{c1})-(\ref{c4}) is examined.
ADOR parameter selections are compared in
FIG. 1 at 4 different times (t=0.3,0.5,0.8,2.0).
These results are all obtained by using 101 grid points
with a time increment of 0.002. 

In FIG. 1e, results for all the four different forms of 
coefficients are plotted  for a detailed comparison. 
For the region away from the shock front, there is essentially 
no difference among four different forms.
Although, results computed from different forms 
are slightly different at the shock front,
four different choices of 
the diffusion coefficients have similar behavior 
and all are capable of correctly simulating Burgers' shock.

The super diffusion coefficient chosen in Eq. (\ref{c3})
does not work very well as the shock front is 
distorted (see  FIG. 1f). However, we have found that 
a combination of 
anisotropic diffusion and super diffusion does work 
better than the choice of a single term 
(results not shown). Similar results were also
obtained by using 50 grid points.
It is mentioned that the tests on the use of 
the edge enhancing functional
does not result in a stable solution with the 
mesh system and time increment, chosen 
the present study.

\subsection{Burgers' equation in two dimensions}

Let us consider Burgers' equation of the form
\begin{eqnarray}\label{BUG1}
&&u_t+uu_x+vu_y={1\over {\rm Re}}(u_{xx}+u_{yy})\\
&&v_t+uv_x+vv_y={1\over {\rm Re}}(v_{xx}+v_{yy})
\end{eqnarray}
in a square $[0,1]\times [0,1]$ with
the initial values
\begin{eqnarray}\label{BUG2}
&&u(x,y,0)=\sin(\pi x)\sin(\pi y)\\
&&v(x,y,0)=[\sin(\pi x)+\sin(2\pi x)]+[\sin(\pi y)+\sin(2\pi y)]
\end{eqnarray}
and boundary conditions
\begin{eqnarray}\label{BUG3}
&&u(0,y,t)=u(1,y,t)=u(x,0,t)=u(x,1,t)=0\\
&&v(0,y,t)=v(1,y,t)=v(x,0,t)=v(x,1,t)=0.
\end{eqnarray}
This problem has no analytical solution and 
is chosen to demonstrate that the present algorithm can be 
efficient, even with a very coarse mesh.
This case is computed by using an ADOR parameter of 
$\gamma^3_1=0.0006$ with $41^2$ points.
The contours of velocity field components at t=1 are 
plotted in FIG. 2.

\subsection{The incompressible Navier-Stokes equations}

To test the present approach for shock capturing further
we consider the Navier-Stokes 
equation
\begin{eqnarray}\label{NS1}
&&u_t+uu_x+vu_y=-p_x+{1\over {\rm Re}}(u_{xx}+u_{yy})\\
&&v_t+uv_x+vv_y=-p_y+{1\over {\rm Re}}(v_{xx}+v_{yy})
\end{eqnarray}
with equation of continuity 
\beq
u_x+v_y=0,
\eeq
where $(u,v)$ is the 
velocity vector, $p$ is the pressure, 
Re (Re$ > 0$) is the Reynolds number
and Re$ =\infty$ defines the Euler equation. 
The domain of problem is 
a square $[0,2\pi]\times [0,2\pi]$ with
periodic boundary conditions. 
With appropriate initial values, the Euler 
equation can be used to describe 
a flow field of vertically perturbed
horizontal shear 
layers around a jet.
Bell et al\cite{BelCol} studied this case 
by a second order projection method. 
Recently E and Shu\cite{EShu} have employed 
this example to demonstrate the
success of their high order ENO scheme 
for resolving the fine vorticity structure  
of the double shear layers.

{\it\bf Case 1: Analytically solvable initial values}

The Navier-Stokes equation is  analytically solvable
for appropriate initial values
\begin{eqnarray}\label{NS5}
&&u(x,y,0)=-\cos(x)\sin(y) \nonumber\\
&&v(x,y,0)=\sin(x)\cos(y).
\end{eqnarray}
The exact solution for this case is given by\cite{EShu}  
\begin{eqnarray}\label{NS6}
&&u(x,y,t)=-\cos(x)\sin(y)e^{-{2t\over {\rm Re}}} \nonumber\\
&&v(x,y,t)=\sin(x)\cos(y)e^{-{2t\over {\rm Re}}} \nonumber\\
&&p(x,y,t)=-{1\over4}[\cos(2x)+\cos(2y)]e^{-{4t\over {\rm Re}}}.
\end{eqnarray}
This provides a benchmark test for potential
numerical methods in fluid dynamics.
The implicit Euler scheme is used for the time 
integration and the DSC algorithm is utilized
for the spatial discretization. 
The accuracy and reliability of this combination 
was previously tested\cite{weicmame20} for 
this problem using a standard LU decomposition 
algorithm for solving linear algebraic equations. 
Here, DSC-ADI algorithm\cite{zhaoADI} is used. 
We choose a grid of $32^2$ for the present calculation with 
a time increment of 0.001.
The DSC-ADI results are summarized in TABLE II.
Note that, for the inviscid case 
(Re=$\infty$) the present result is 
accurate to the machine precision.
It is evident that the accuracy of the DSC 
approach is extremely high, particularly when 
the Reynolds numbers are very large.

{\it\bf Case 2: The Euler equation}

We now test our anisotropic diffusion approach 
for the Euler equation (Re$ =\infty$) 
with sharply varying initial values.
This example is chosen to 
illustrate the ability of the present approach, 
for providing very fine resolution with a 
relatively coarse grid. 
The initial values are that of a jet in a doubly 
periodic geometry
\begin{eqnarray}\label{NS7}
&&u(x,y,0)=\left\{\begin{array}{ll}\mbox{$
\tanh\left({2y-\pi\over 2\rho}\right),$} & \mbox{if $y\leq\pi$}\\
\mbox{$\tanh\left({3\pi-2y\over2\rho}\right),$} & \mbox{if $y >\pi$}
\end{array}
\right\} \nonumber\\
&&v(x,y,0)=\delta\sin(x),
\end{eqnarray}
where $\delta=0.05$ is used for the convenience of 
comparison with the previous study\cite{EShu}.
This initial value describes the flow field consisting of 
horizontal shear layers of finite thickness, perturbed by a small 
amplitude vertical velocity, making up the boundaries of the jet. 
However, this problem is not analytically solvable.
A pioneer work in this study was given by 
Bell et al\cite{BelCol}, in which 
they utilized a second-order Godunov scheme in association with 
a projection approach for divergence-free velocity fields
with general boundary condition. With a periodic boundary condition,
E and Shu have shown that their high-order ENO scheme\cite{EShu}
performs well for this problem.

We first consider the parameter  $\rho=\pi/15$, 
a case studied by Bell et al\cite{BelCol} using 
a projection method with three sets of grids
(128$^2$, 256$^2$ and  512$^2$).
E and Shu\cite{EShu} computed this case by 
using both  spectral collocation code 
with 512$^2$ points and their high order ENO scheme
with 64$^2$ and  128$^2$
points. The spectral collocation code
produced an oscillatory 
solution at $t=10$ (see FIG. 1 of Ref. \cite{EShu}), 
while the high order ENO scheme
produced a defect at $t=6$ as the channels 
connecting the vorticity centers are slightly
distorted (see FIG. 2 of Ref. \cite{EShu}).
In the present simulation, we
choose a $64^2$ grid for the 
computational domain with a time increment of 0.002.
The prescription in Eq. (\ref{c4}) is used with the 
ADOR parameter of
$\gamma^3_1=0.0006$. 
The results at different times (t=4,6,8,10)
are plotted in FIG 3. 
It is seen that our solutions are smooth (some non-smooth
features in the contour plot is due to the fact 
that the grid is very coarse) and stable for this case. 
In particular, no distortion is found in vorticity contours at t=6.
For early times, present results compare extremely 
well with those of the spectral collocation code computed 
with 512$^2$ points. There are no spurious
numerical oscillations during the entire process.

Finally, we perform simulations by the  
 present approach using 
the discontinuous initial data ($\rho\rightarrow0$)
as in \cite{EShu}. 
The evolution under the Euler equation 
leads to a periodic array of large vortices, 
with the shear layers between the rolls being 
thinned by the large straining field.   
It is known that, for the present choice of parameters
the solution quickly develops into a roll-up process with smaller 
and smaller scales, and the resolution is lost eventually 
with a fixed grid for local methods.
We compute this case by using 
$128^2$ grid points as in \cite{EShu}. A time increment  
of 0.002 is used.  The ADOR parameter 
is chosen as $\gamma^3_1=0.0006$. 
The present results at different times (t=4,6,8,10,12,14)
are plotted in FIG 4.

\section{Conclusions}

Connection is made between digital image processing 
and computational fluid dynamics. The evolution of 
an image surface under a partial differential operator 
can be viewed as a form of image processing. 
Computationally, numerical shock capturing can be 
formulated on the lines of iterative edge-detection.
Hence, techniques developed in the computational fluid 
dynamics can be used for image processing and 
vice-versa. This paper introduces 
the method of anisotropic diffusion
oscillation reduction (ADOR), an approach which has its roots
in image processing, for shock wave computations. 

In fact, the ADOR method is much similar to the 
artificial viscosity algorithm. 
Physical origins and mathematical 
properties of the artificial viscosity are  
discussed from the kinetic theory point of view.
The form of pressure tensor is derived from the 
first principles of quantum mechanics. 
Quantum kinetic theory is utilized 
to arrive at macroscopic transport equations from the 
microscopic quantum theory.  Macroscopic symmetry is used
to simplify the phenomenological pressure tensor expressions. 
The latter provides the basis for the design of artificial 
viscosity. The original von Neumann-Richtmyer form
of artificial viscosity fits into a special case of 
the present generalizations.

The anisotropic diffusion, which is essentially
an image processing technique, 
is modified for shock capturing. The technique 
preserves image edges by introducing little diffusion, 
where the image gradient is large, while it provides
a substantial diffusion coefficient at smooth parts 
of the image. In the present shock capturing algorithm,
the edge-sensitive diffusion is introduced at a 
shock front so that large oscillations can be 
efficiently eliminated. An edge enhancing functional 
and edge-detected super diffusion operators are 
proposed for shock capturing. These terms are  
introduced from the general structure of conservation 
laws. In fact, kinetic theory analysis and anisotropic 
diffusion argument lead to a number of similar 
expressions for shock wave treatments. 
Hence, the acronym ADOR is referred for both 
approaches.

The reliability and robustness of the ADOR method is
explored in association with the discrete singular 
convolution (DSC) algorithm\cite{weijcp99,weiphysica20,weijpb20,weijpa20}. 
The DSC algorithm was previously validated by handling many 
linear and nonlinear problems and is a potential algorithm 
for the computer realization of some singular convolutions.
A few detailed prescriptions for the anisotropic 
diffusion functional are considered in this work.
The coefficients in these prescriptions
are chosen to be of lower order in gradient so that 
a coarse grid can be used. 
A number of standard test examples, including 
(inviscid) Burgers' equation in one and two spatial 
dimensions, the incompressible Navier-Stokes and the Euler 
equations,  are employed for the present test computations. 

For  Burgers' equation, we have examined the accuracy for  
spatial and temporal discretizations. A few ADOR coefficients
given by Eqs. (\ref{c1})-(\ref{c4}) are tested. 
The anisotropic diffusion term is found to be very robust 
in association with many edge-detected coefficients,
while the super diffusion operator produces an over shot 
at Burgers' shock front. The edge enhancing 
functional apparently does not perform well for this problem.
The ADOR approach is found to be very stable for the 
2D inviscid Burgers' equation even with a very coarse grid.

For the incompressible Navier-Stokes equation, the ADOR 
approach is used in association with a DSC-ADI algorithm. 
The accuracy of the algorithm was tested by an analytically
solvable case and the machine precision is found at 
the inviscid limit of the Navier-Stokes equation. 
The algorithm was then applied to the study of  
doubly periodic shear layers, which was chosen to 
demonstrate the capability of the ADOR method for more 
difficult problems. These results indicate that 
the proposed method has the potential to capture shock waves. 
Work on the implementation of the 
present approach to more complex fluid flow systems, 
including general boundary conditions and for 
compressible flow is under progress.

\vspace*{1cm}

\centerline{\bf Acknowledgment}
 
{This work was supported in part by the 
National University of Singapore.
The author thanks Professors Tony Chan and C.-W. Shu 
for useful discussions
about the ENO and WENO schemes.}

\vspace*{1cm}

      





\begin{table}
\caption{ $L_1$  and $L_\infty$  
errors of the DSC solutions for 
           Burgers' equation}
\begin{center}
\begin{tabular}{c|c|c} 
Time   & $L_{1}$ & $L_{\infty}$ \\ \hline
  0.4  &    2.1(-4) &  2.8(-3)\\
  0.8  &    2.8(-4) &  4.4(-3)\\
   1.2 &    4.2(-5) &  6.8(-4)\\
   1.6 &    5.8(-6) &  9.4(-5)\\
   2.0 &    9.4(-7) &  1.2(-5)\\
   3.0 &    3.8(-8) &  4.0(-7)\\
   5.0 &    2.1(-10)&  2.2(-9)\\
   10  &    1.5(-11)&  4.0(-11)\\
   20  &    3.2(-12)&  5.3(-12)\\
   30  &    1.1(-12)&  1.8(-12)\\
   40  &    4.1(-13)&  6.4(-13)
\end{tabular}
\end{center}               
\end{table}

\begin{table}
\caption{ $L_2$  errors of the DSC-ADI solutions for 
            the 2D Navier-Stokes equation computed
           with 32 grid point in each dimension}
\begin{center}
\begin{tabular}{c|c|c|c|c|c|c|c|c}
 Re  & \multicolumn{2}{c|}{t=1}& \multicolumn{2}{c|}{t=2}
     & \multicolumn{2}{c|}{t=3}& \multicolumn{2}{c}{t=4}
\\ \cline{2-9}
     &  $u$ & $v$ &  $u$ & $v$ &  $u$ & $v$ &  $u$ & $v$ 
 \\ \hline 
      $20$  &
  6.1(-11)&
  6.1(-11)& 
  1.1(-10)$^a$ &  
  1.1(-10)& 
  1.5(-10)&
  1.5(-10)&
  1.8(-10)&
  1.8(-10) 
\\
      $10^2$  &
  1.2(-12)& 
  1.2(-12)&  
  2.3(-12)& 
  2.3(-12)&  
  3.4(-12)& 
  3.4(-12)&  
  4.4(-12)& 
  4.4(-12)
\\
      $10^3$  & 
  9.5(-13)& 
  9.3(-13)& 
  1.9(-12)& 
  1.8(-12)& 
  2.8(-12)& 
  2.8(-12)& 
  3.8(-12)& 
  3.7(-12) 
\\
      $10^4$  & 
  8.4(-13)& 
  8.3(-13)& 
  1.7(-12)& 
  1.7(-12)& 
  2.6(-12)& 
  2.6(-12)& 
  3.5(-12)& 
  3.5(-12) 
\\
      $10^5$  & 
  4.6(-13)& 
  4.6(-13)&
  9.7(-13)& 
  9.7(-13)&
  1.4(-12)& 
  1.4(-12)&
  1.9(-12)& 
  1.9(-12)
\\
      $10^6$  & 
  6.3(-13)& 
  6.3(-13)& 
  1.2(-12)& 
  1.2(-12)& 
  1.7(-12)& 
  1.9(-12)& 
  2.6(-12)& 
  2.6(-12)
\\
      $\infty$& 
  1.5(-15) & 
  1.5(-15) &
  1.6(-15) $^b$ & 
  1.6(-15) & 
  1.6(-15) & 
  1.6(-15) & 
  1.7(-15) & 
  1.7(-15)  
\\ 
\end{tabular}
\end{center}
$^a$: $L_2=9.1(-04)$ and $^b$: $L_2=4.9(-04)$ 
   were reported in Ref. \cite{EShu}.   
\end{table}

\newpage

\centerline{\bf Figure Captions}

{\bf FIG. 1.} A comparison of various ADOR coefficients
for the numerical solution of inviscid Burgers' equation. 
(a) $d^1_1=0.0009,d^1_2=0,e^1=0$;
(b) $\gamma^1_1=0.0022,\gamma^1_2=0$;
(c) $\gamma^2_1=0.0018,\gamma^2_2=0$;
(d) $\gamma^3_1=0.0015,\gamma^3_2=0$;
(e) Comparison of (a), (b),(c) and (d);
(f) $\gamma^3_1=0,\gamma^3_2=-98\times 10^{-8}$.

\vskip 20pt

{\bf FIG. 2.} 
The ADOR solution for the 2D inviscid Burgers' equation
with 41$^2$ points, $\gamma^4_1=0.02, \gamma^4_2=0$.  
Left: the $u$ field; right: the $v$ field.

\vskip 20pt

{\bf FIG. 3.} 
The vorticity contours for 
the 2D Euler equation  by the ADOR method with 
64$^2$ points, $\gamma^5_1=0.0006$. 
Up left:   t=4;
up right:  t=6;
low left:  t=8;
low right: t=10.

\vskip 20pt

{\bf FIG. 4.} 
The vorticity contour for the 2D Euler equation, 
with discontinuous initial data,  by the ADOR 
method with 128$^2$ points, $\gamma^5_1=0.0006$. 
Up left:   t=4;
up right:  t=6;
middle left:  t=8;
middle right: t=10;
low left:  t=12;
low right: t=14.

\newpage

\begin{center}
\resizebox{17cm}{!}{\includegraphics{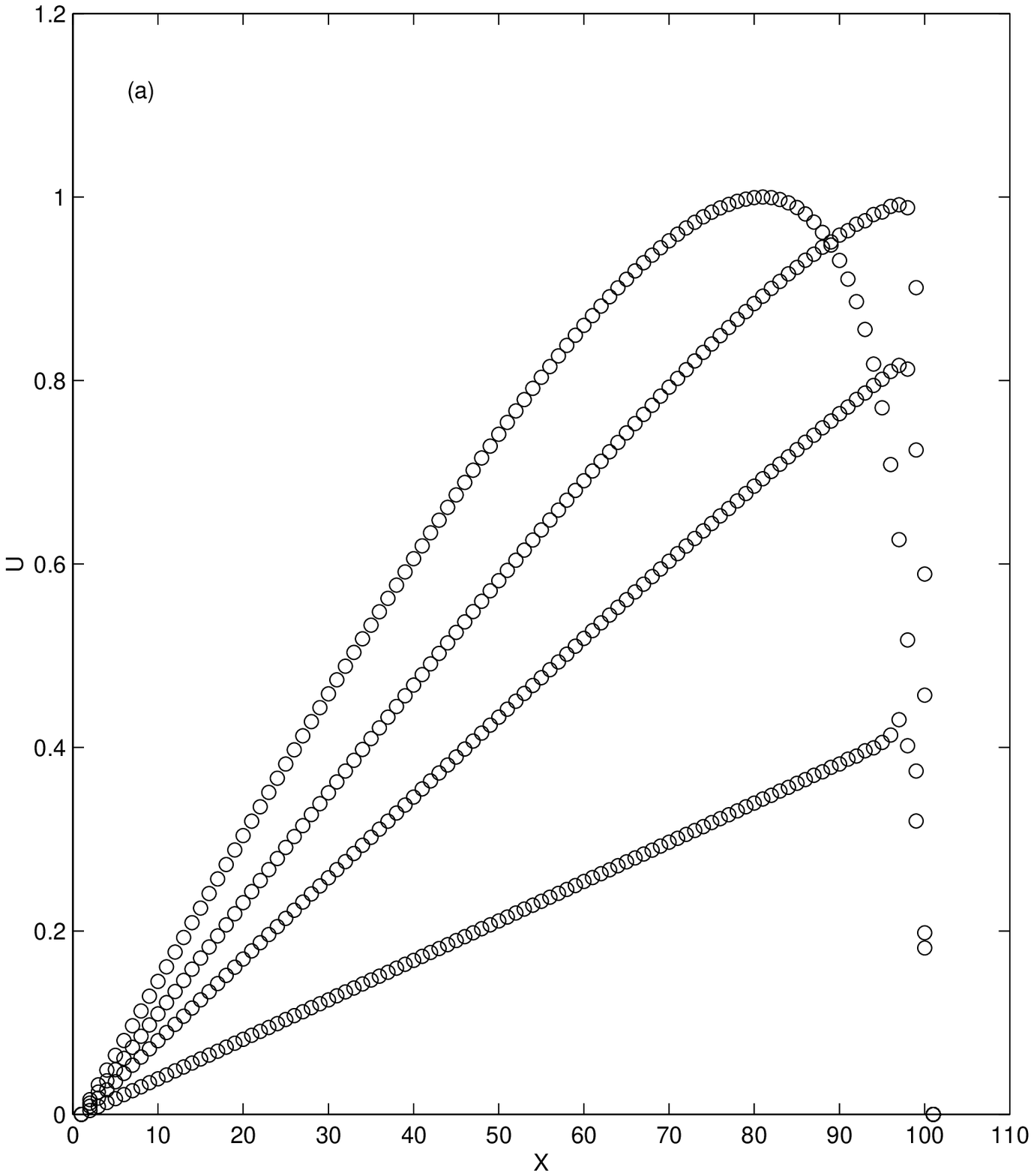}}

\resizebox{17cm}{!}{\includegraphics{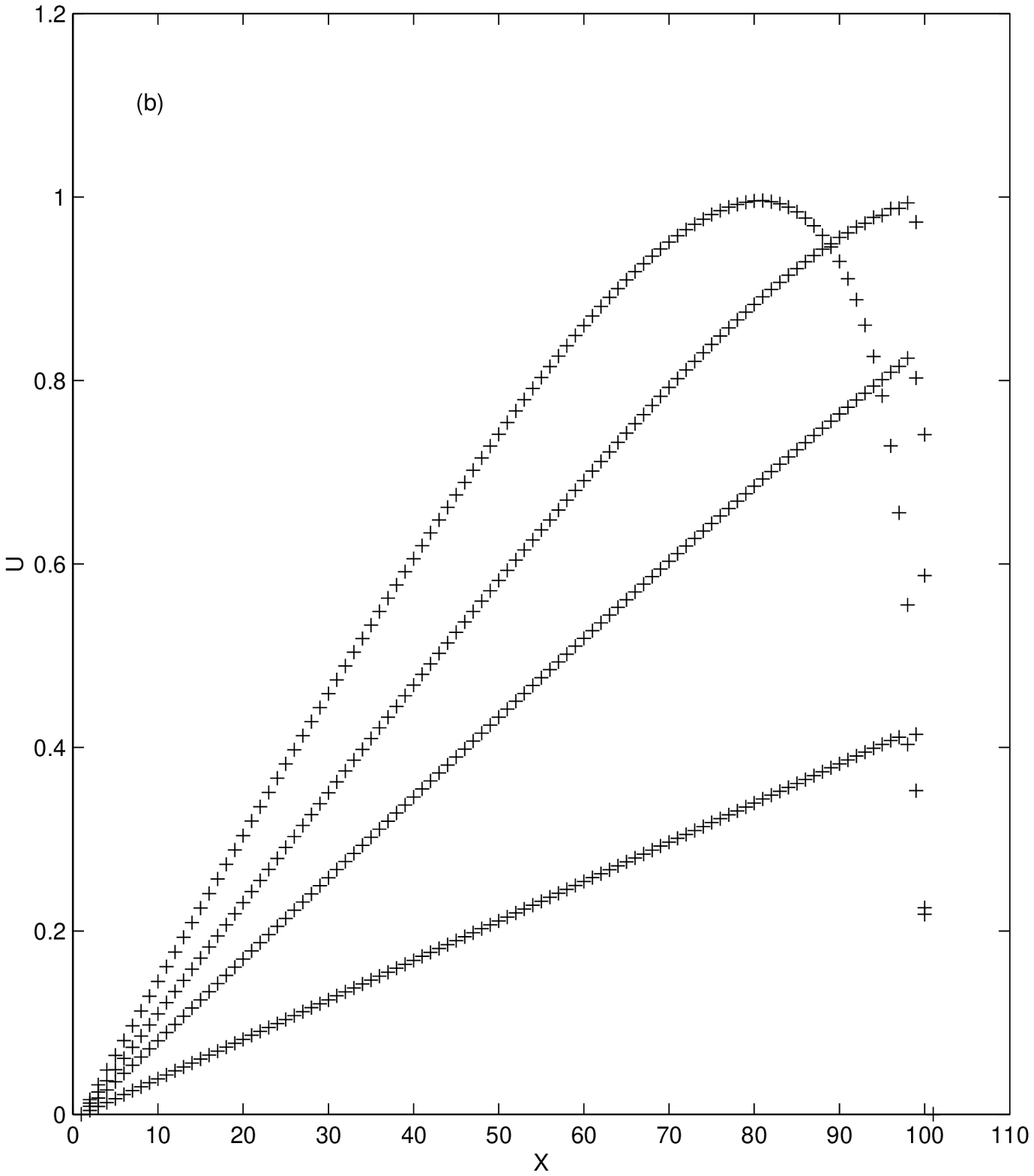}}

\resizebox{17cm}{!}{\includegraphics{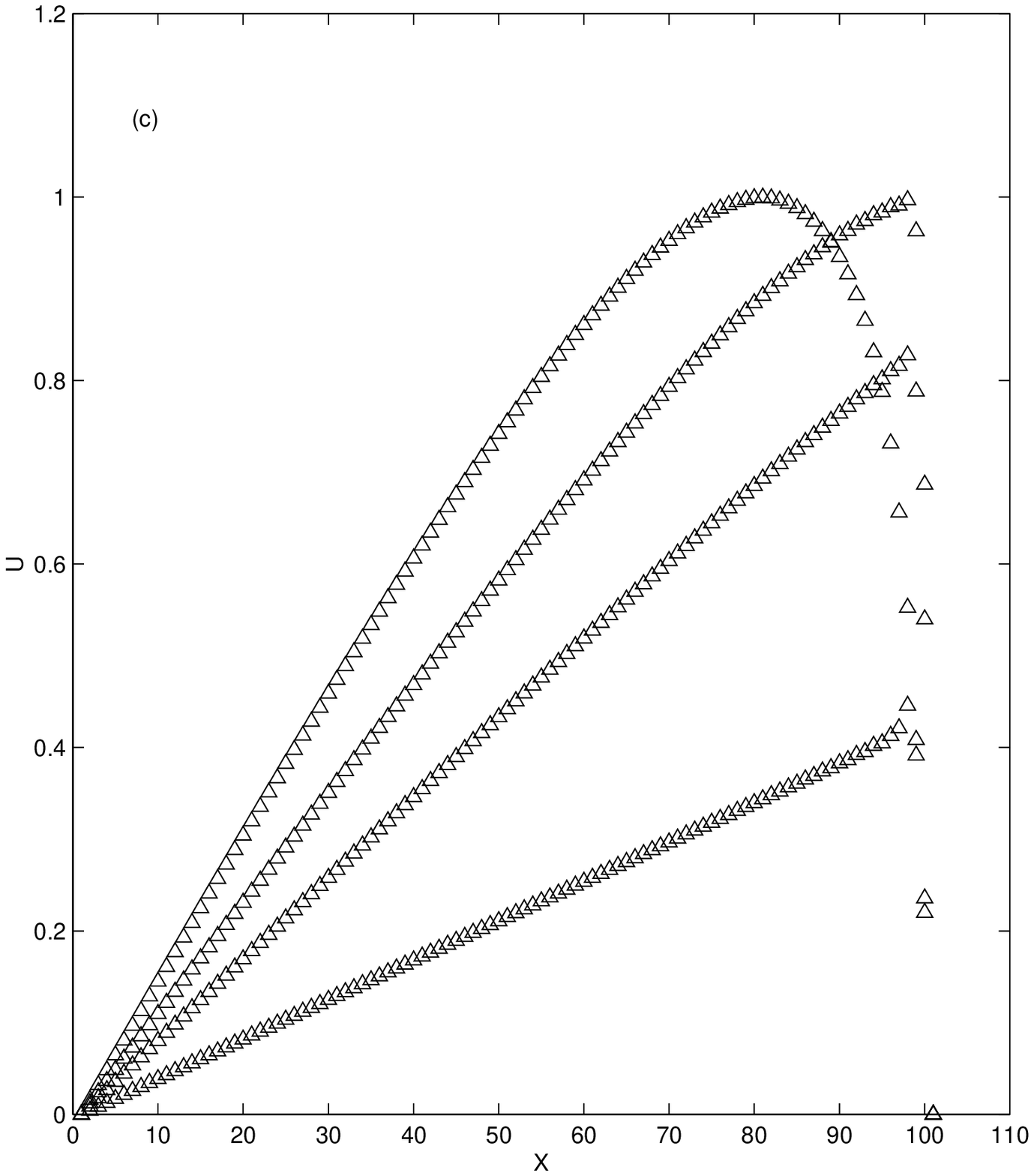}}

\resizebox{17cm}{!}{\includegraphics{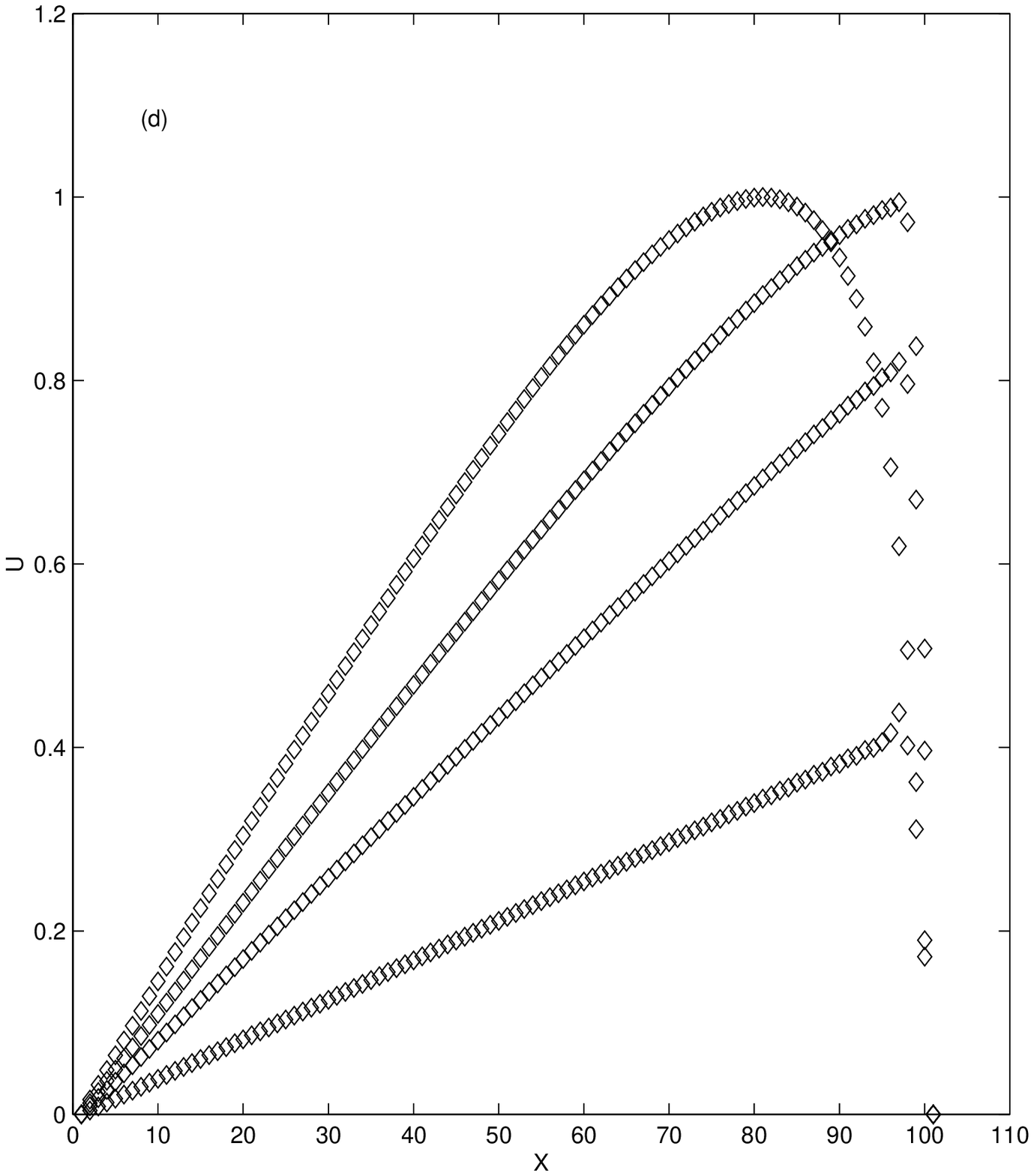}}

\resizebox{17cm}{!}{\includegraphics{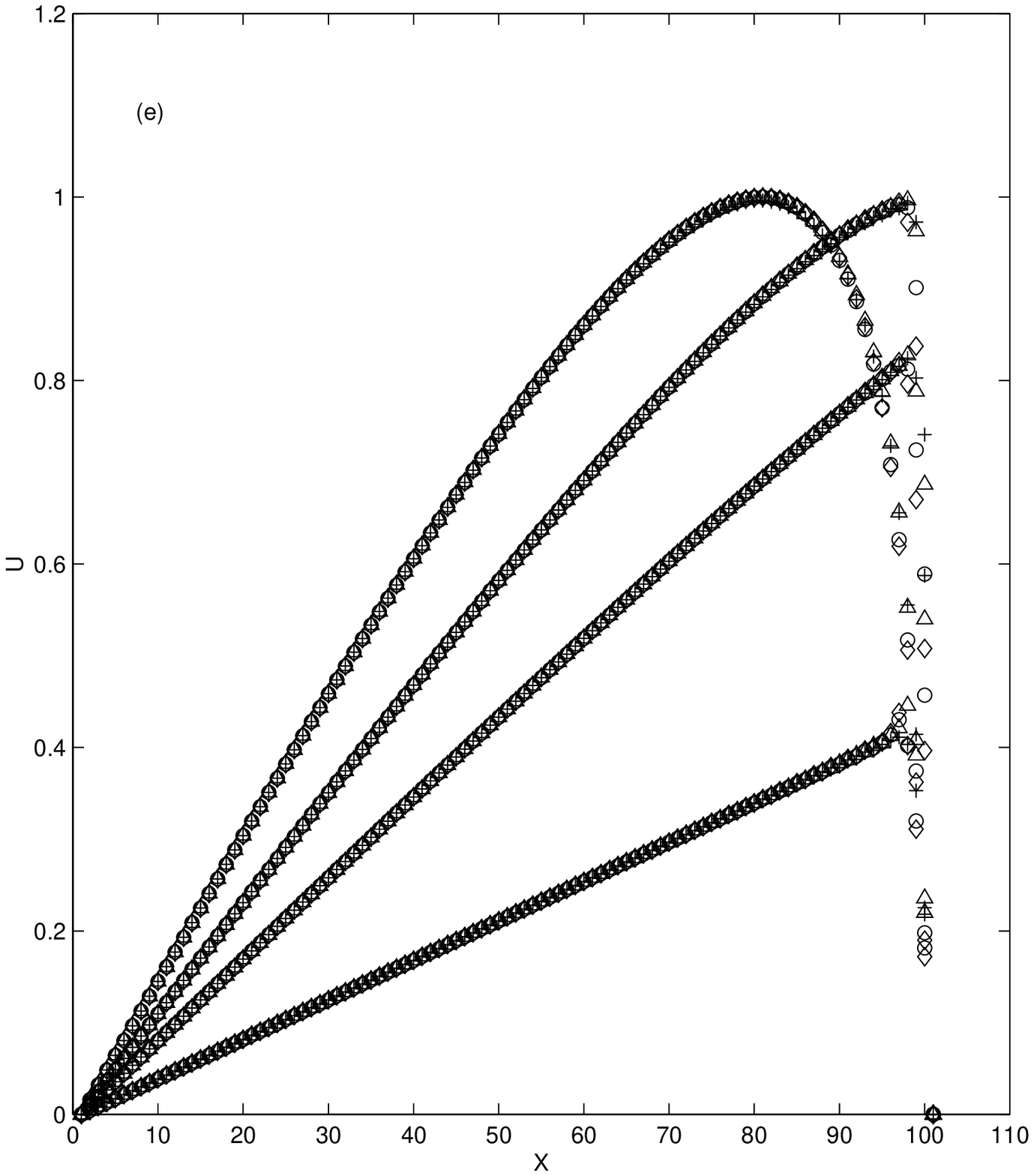}}

\resizebox{17cm}{!}{\includegraphics{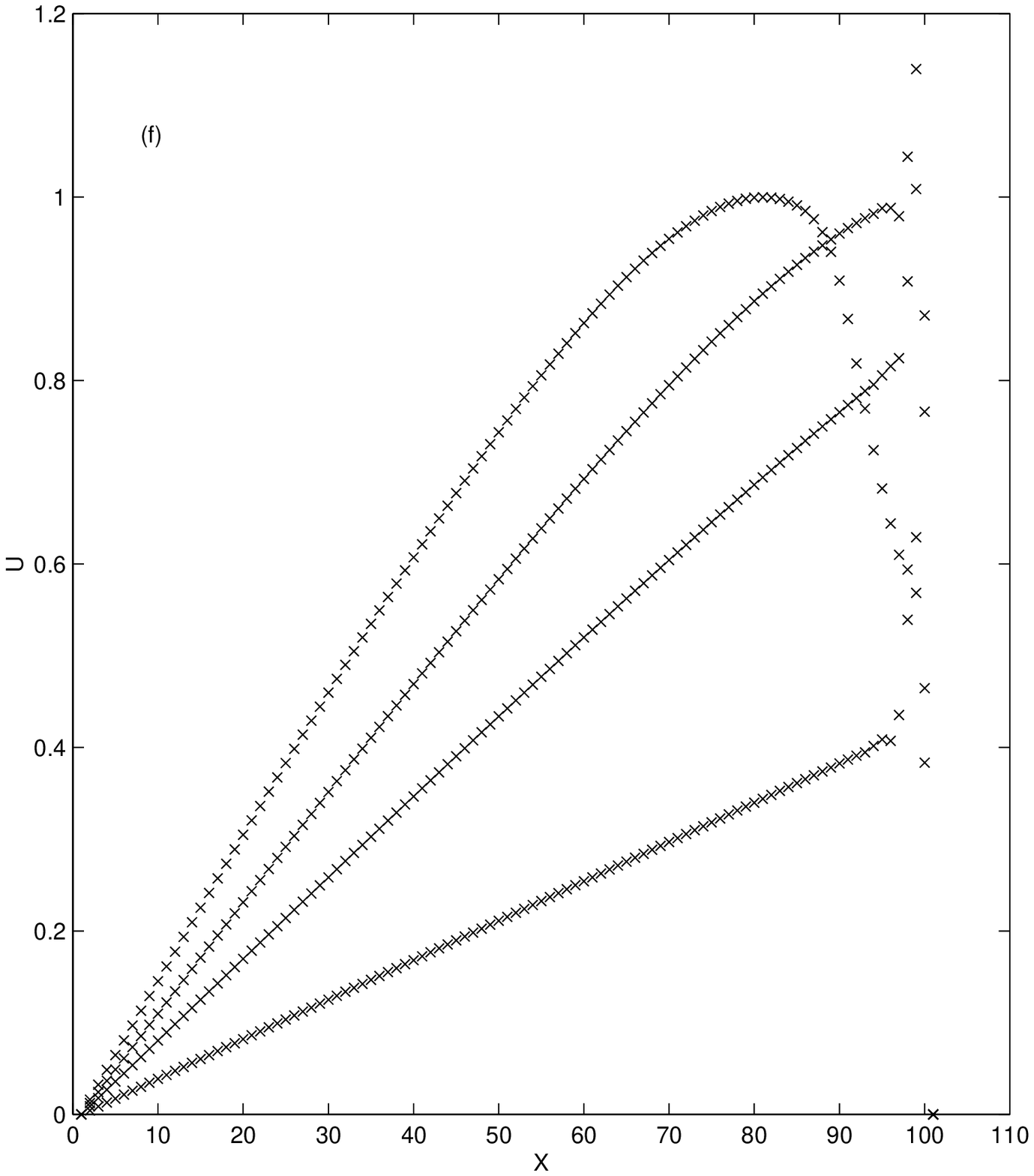}}
\end{center}
{\bf FIG. 1.} A comparison of various ADOR coefficients
for the numerical solution of inviscid Burgers' equation. 
(a) $d^1_1=0.0009,d^1_2=0,e^1=0$;
(b) $\gamma^1_1=0.0022,\gamma^1_2=0$;
(c) $\gamma^2_1=0.0018,\gamma^2_2=0$;
(d) $\gamma^3_1=0.0015,\gamma^3_2=0$;
(e) Comparison of (a), (b),(c) and (d);
(f) $\gamma^3_1=0,\gamma^3_2=-98\times 10^{-8}$.

\newpage

\begin{center}
\resizebox{17cm}{!}{\includegraphics{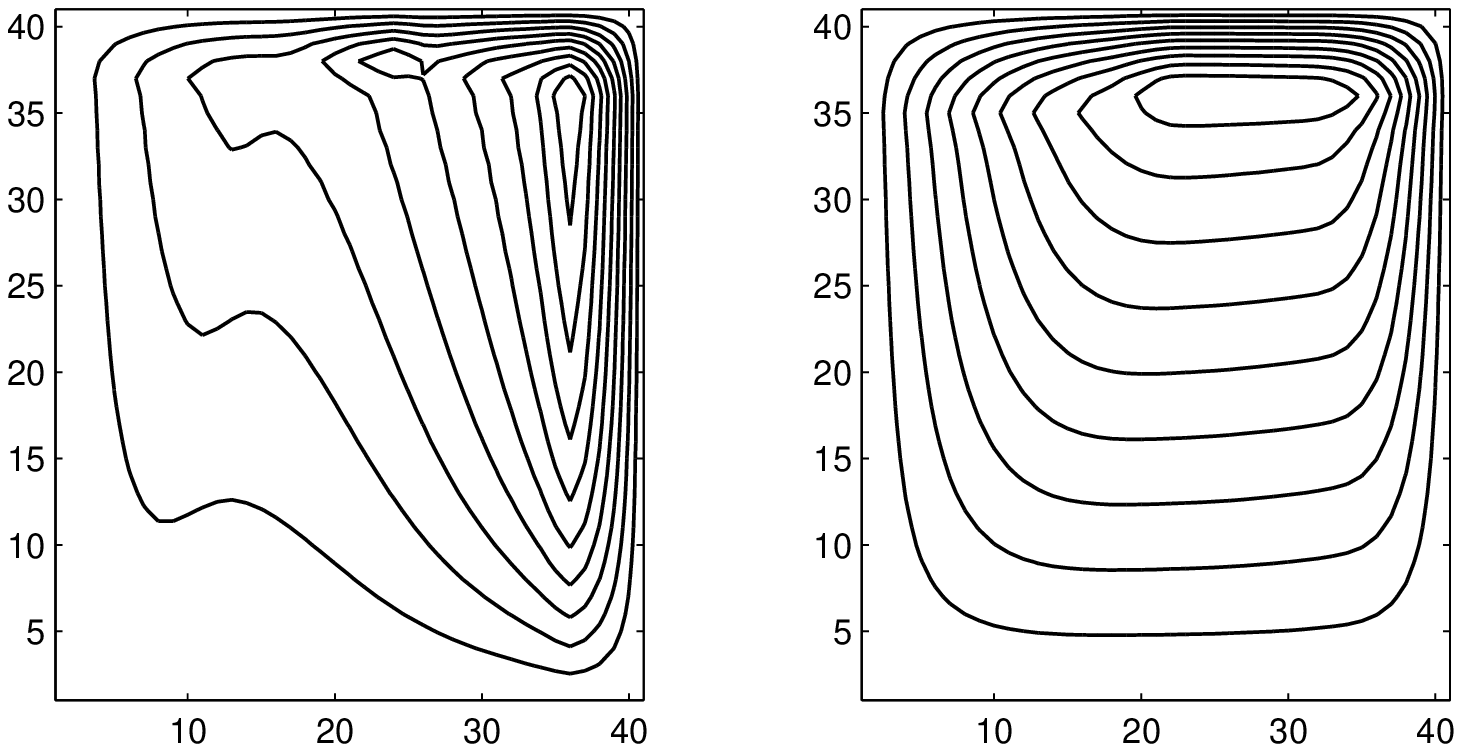}}
\end{center}
{\bf FIG. 2.} 
The ADOR solution for the 2D inviscid Burgers' equation
with 41$^2$ points, $\gamma^4_1=0.02, \gamma^4_2=0$.  
Left: the $u$ field; right: the $v$ field. 

\newpage

\begin{center}
\resizebox{17cm}{!}{\includegraphics{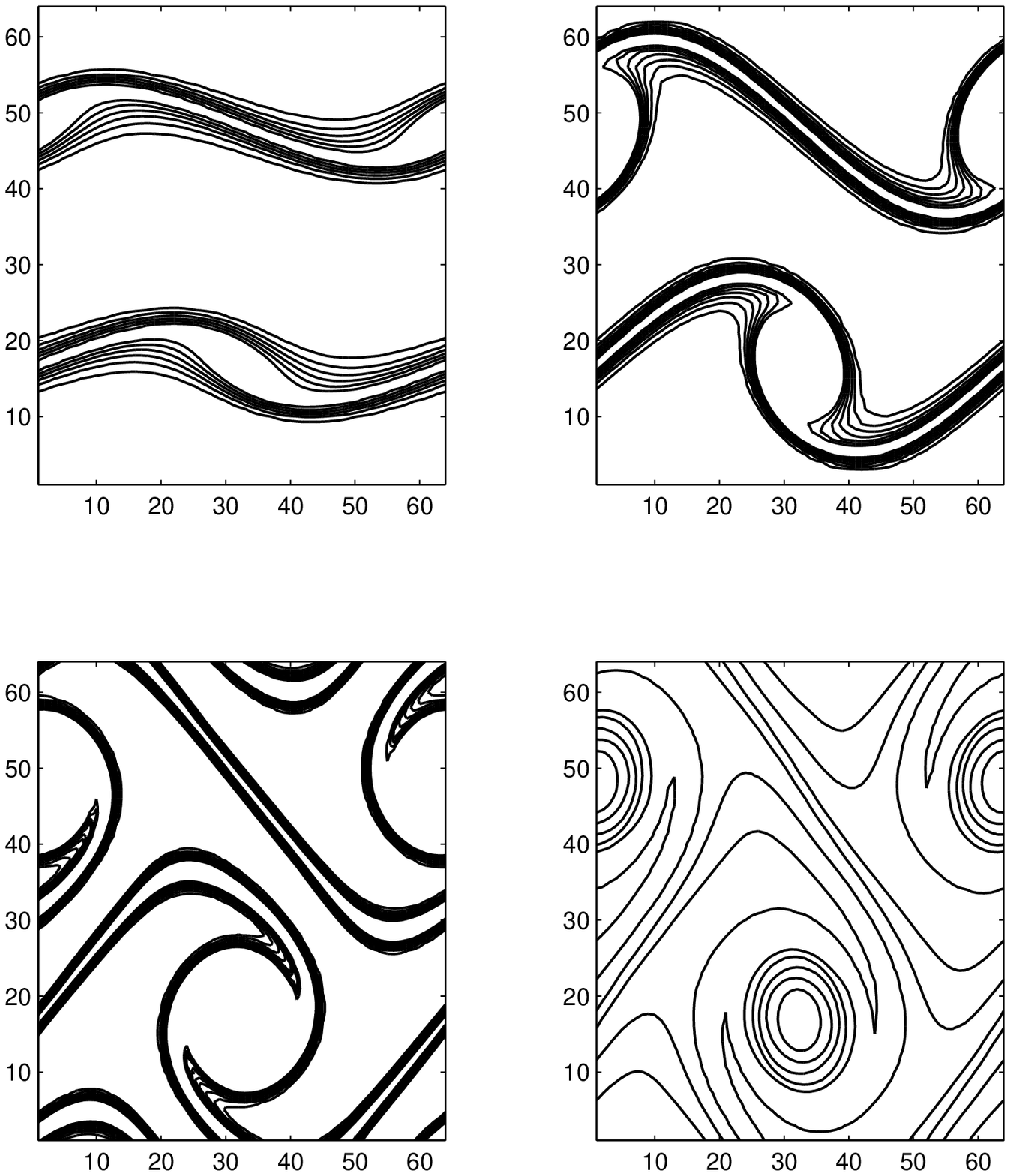}}
\end{center}
{\bf FIG. 3.} 
The vorticity contours for 
the 2D Euler equation  by the ADOR method with 
64$^2$ points, $\gamma^5_1=0.0006$. 
Up left:   t=4;
up right:  t=6;
low left:  t=8;
low right: t=10.

\newpage

\begin{center}
\resizebox{17cm}{!}{\includegraphics{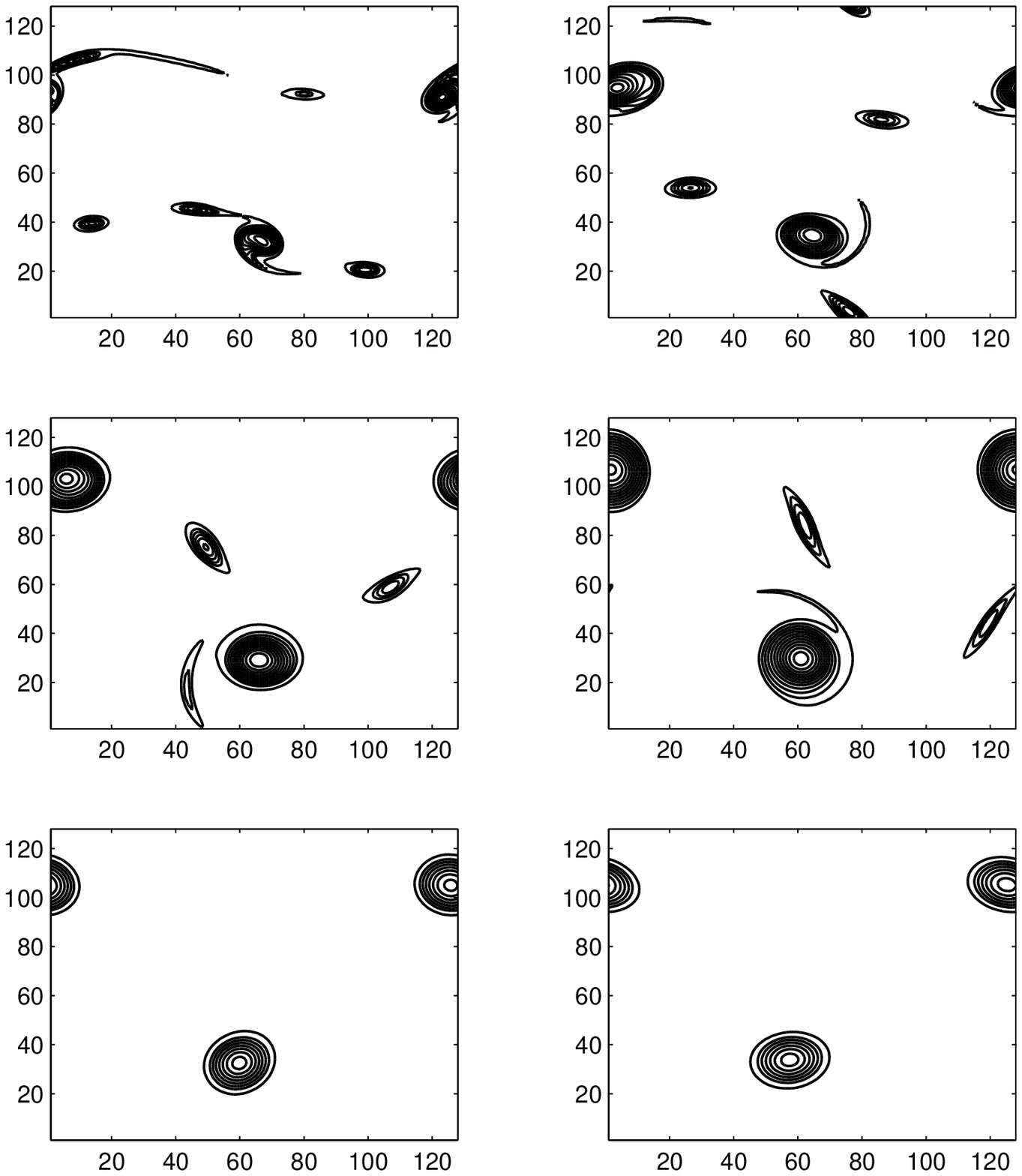}}
\end{center}
{\bf FIG. 4.} 
The vorticity contour for the 2D Euler equation, 
with discontinuous initial data,  by the ADOR 
method with 128$^2$ points, $\gamma^5_1=0.0006$. 
Up left:   t=4;
up right:  t=6;
middle left:  t=8;
middle right: t=10;
low left:  t=12;
low right: t=14.

\end{document}